# Covering Constants for Metric Projection Operator with Applications to Stochastic Fixed-Point Problems

Jinlu Li

Department of Mathematics
Shawnee State University
Portsmouth, Ohio 45662
USA

**Abstract**

In this paper, we use the Mordukhovich derivatives to precisely find the covering constants for the metric projection operator onto nonempty closed and convex subsets in uniformly convex and uniformly smooth Banach spaces. We consider three cases of the subsets: closed balls in uniformly convex and uniformly smooth Banach spaces, closed and convex cylinders in $l_p$, and the positive cone in $L_p$, for some $p$ with $1 < p < \infty$. By using Theorem 3.1 in [2] and as applications of covering constants obtained in this paper, we prove the solvability of some stochastic fixed-point problems. We also provide three examples with specific solutions of stochastic fixed-point problems.

## 1. Introduction

Let $X$ and $Y$ be Banach spaces and let $U$ and $V$ be nonempty subsets in $X$ and $Y$, respectively. Let $\Phi: X \rightrightarrows Y$ be a multifunction (a set-valued mapping). The graph of $\Phi$ is denoted by gph $\Phi$, which is a subset of $X \times Y$ defined by

$$\text{gph } \Phi \coloneqq \{(x, y) \in X \times Y : y \in \Phi(x)\}.$$

We say that $\Phi$ enjoys the covering property with modulus $\alpha > 0$ (or it has the $\alpha$-covering property) on $U$ relative to $V$ if (see (2.1) in [2])

$$\Phi(x) \cap V + \alpha r \mathbb{B}_Y \subset \Phi(x + r\mathbb{B}_X), \text{ whenever } x + r\mathbb{B}_X \subset U, \text{ as } r > 0. \qquad (1.1)$$

Here, $\mathbb{B}_X$ and $\mathbb{B}_Y$ are the closed unit balls in $X$ and $Y$, respectively. Notice that, in [2], the covering property of set-valued mappings is defined in normed spaces. The covering property for set-valued mappings is a very important concept in variational analysis and optimization (see (see [2−10, 21−22]).

By the covering criterion (see Theorem 2.1 in [2]), there is a very close connection between the covering property and covering constant for set-valued mappings. The covering constant is defined with respect to the Mordukhovich derivatives of the considered mappings. So, we need to review the concept of Mordukhovich derivative of the set-valued mapping $\Phi$.

Let $(\bar{x}, \bar{y}) \in \text{gph } \Phi$, the Mordukhovich derivative, or Mordukhovich coderivative, or coderivative of $\Phi$ at $(\bar{x}, \bar{y})$, denoted by $\widehat{D}^*\Phi(\bar{x}, \bar{y})$, is a set-valued mapping $\widehat{D}^*\Phi(\bar{x}, \bar{y}): Y^* \rightrightarrows X^*$ that is defined, for any $y^* \in Y^*$, by (see Definitions 1.13 and 1.32 in [17])

$$\widehat{D}^*\Phi(\bar{x}, \bar{y})(y^*) = \left\{ x^* \in X^*: \limsup_{\substack{(u,v) \to (\bar{x},\bar{y}) \\ (u,v) \in \text{gph } \Phi}} \frac{\langle x^*, u-\bar{x}\rangle_X - \langle y^*, v-\bar{y}\rangle_Y}{\|u-\bar{x}\|_X + \|v-\bar{y}\|_Y} \leq 0 \right\}. \tag{1.2}$$

Here, $\|\cdot\|_X$ and $\|\cdot\|_Y$ denote the norms in $X$ and $Y$, respectively. $\langle \cdot, \cdot \rangle_X$ denotes the natural paring between $X^*$ and $X$; $\langle \cdot, \cdot \rangle_Y$ denotes the natural paring between $Y^*$ and $Y$.

Mordukhovich derivative of set-valued mappings has laid the foundation of generalized differentiation in set-valued analysis, which has been widely applied to operator theory, optimization theory, approximation theory, control theory, equilibrium theory, and so forth (see [16−19]). One of the most important applications of Mordukhovich derivatives of set-valued mappings is to to define the covering constants for the considered set-valued mappings. The covering constant for $\Phi$ at point $(\bar{x}, \bar{y}) \in \text{gph } \Phi$ is defined by (see (2.6) in [2])

$$\hat{\alpha}(\Phi, \bar{x}, \bar{y}) := \sup_{\eta > 0} \inf\{\|z^*\|_{X^*}: z^* \in \widehat{D}^*\Phi(x,y)(w^*), x \in \mathbb{B}_X(\bar{x}, \eta), y \in \Phi(x) \cap \mathbb{B}_Y(\bar{y}, \eta), \|w^*\|_{Y^*} = 1\}. \tag{1.3}$$

Here, $\|\cdot\|_{X^*}$ and $\|\cdot\|_{Y^*}$ denote the norms in $X^*$ and $Y^*$, respectively. $\mathbb{B}_X(\bar{x}, \eta)$ is the closed ball in $X$ centered at $\bar{x}$ with radius $\eta$, and $\mathbb{B}_Y(\bar{y}, \eta)$ is the closed ball in $Y$ centered at $\bar{y}$ with radius $\eta$. Now we review the covering criterion.

**Theorem 2.1 in [2].** *Let $\Phi: X \rightrightarrows Y$ be a set-valued mapping between normed spaces, let $(\bar{x}, \bar{y}) \in \text{gph } \Phi$, and let $\hat{\alpha}(\Phi, \bar{x}, \bar{y})$ be taken from* (1.3) *(it is (2.6) in [2]). Then*:

(i) *The $\alpha$-covering property of $\Phi$ around $(\bar{x}, \bar{y})$ implies that $\alpha \leq \hat{\alpha}(\Phi, \bar{x}, \bar{y})$.*

(ii) *Assume that $X$ and $Y$ are Asplund spaces. If $0 < \alpha < \hat{\alpha}(\Phi, \bar{x}, \bar{y})$, then $\Phi$ has the $\alpha$-covering property around $(\bar{x}, \bar{y})$ provided that its graph is closed around this point.*

In the fields of variational analysis and optimization theory related to set-valued mappings, the concept of covering constant (or the covering property) of the considered mappings plays a very important role in the proof of the existence of solutions of parameterized generalized equations, implicit function and fixed-point theorems, optimal value functions in parametric optimization, and so forth (see [2−10]]. For example, we list the main theorem in [2], which proves the existence of parameterized coincidence points for some set-valued mappings.

**Theorem 3.1 in [2].** *Let the spaces $X$ and $Y$ in be Asplund and let $P$ be a topological space. Let the following conditions be satisfied*:

(A1) *The multifunction $F$ is closed around $(\bar{x}, \bar{y})$.*
(A2) *There are neighborhoods $U \subset X$ of $\bar{x}$, $V \subset Y$ of $\bar{y}$, and $O$ of $\bar{p} \in P$ as well as a number $l \geq 0$ such that the multifunction $G(\cdot, p)$ is Lipschitz-like on $U$ relative to $V$ for each $p \in O$ with*

*the uniform modulus l, while the multifunction $p \to G(\bar{x}, p)$ is lower/inner semicontinuous at $\bar{p}$.*

(A3) *The Lipschitzian modulus l of $G(\cdot, p)$ is chosen as $l < \hat{\alpha}(F, \bar{x}, \bar{y})$, where $\hat{\alpha}(F, \bar{x}, \bar{y})$ is the covering constant of F around $(\bar{x}, \bar{y})$ taken from (2.6) in [2].*

*Then for each $\alpha > 0$ with $l < \alpha < \hat{\alpha}(F, \bar{x}, \bar{y})$, there exist a neighborhood $W \subset P$ of $\bar{p}$ and a single-valued mapping $\sigma: W \to X$ such that whenever $p \in W$ we have*

$$F(\sigma(p)) \cap G(\sigma(p), p) \neq \emptyset \quad \text{and} \quad \|\sigma(p) - \bar{x}\|_X \leq \frac{\text{dist}(\bar{y}, G(\bar{x}, p))}{\alpha - l}.$$

This theorem demonstrates the crucial importance of the covering constant for a considered mapping to guarantee the existence of solutions for some parameterized generalized equations. This leads us to find the exact covering constants for some mappings, which are important in set-valued and variational analysis.

In operator theory, the metric projection operator is one of the most useful and important operators. In this paper, we only consider the metric projection operator onto nonempty closed and convex subsets in uniformly convex and uniformly smooth Banach spaces, in which the metric projection operator is a continuous single-valued mapping.

This paper is organized as follows: In section 2, we overview some needed fundamental notations and properties of the metric projection onto nonempty closed and convex subsets in uniformly convex and uniformly smooth Banach spaces and formulate the formula to calculate the covering constant for the metric projection. By using the results from [13−15] about the Mordukhovich derivatives of the metric projection onto closed and convex subsets in uniformly convex and uniformly smooth Banach spaces, we find the exact solutions of the covering constants for the metric projection. More precisely speaking, in section 3, we find precisely the covering constants for the metric projection onto closed balls in uniformly convex and uniformly smooth Banach spaces; in section 4, we find the exact covering constants for the metric projection onto closed and convex cylinders in $l_p$, for $p$ with $1 < p < \infty$; in section 5, we find the exact covering constants for the metric projection onto the positive cone in $L_p$, for $p$ with $1 < p < \infty$. In section 6, by using Theorem 3.1 in [2] and as applications of covering constants obtained in this paper, we prove the solvability of some stochastic fixed-point problems. We also provide three concrete examples with specific solutions of stochastic fixed-point problems. Two of them are deal with single-valued mappings and one of them is for set-valued mappings.

2. **Preliminary**

Let $(X, \|\cdot\|)$ be a real uniformly convex and uniformly smooth Banach space with topological dual space $(X^*, \|\cdot\|_*)$. Since every reflexive Banach space is Asplund (see [2, 20]), so $X$ is Asplund. Let $\langle \cdot, \cdot \rangle$ denote the real canonical pairing between $X^*$ and $X$. Let $\theta$ and $\theta^*$ denote the origins in $X$ and $X^*$, respectively. The identity mappings on $X$ and $X^*$ are respectively denoted by $I_X$ and $I_{X^*}$. Let $J: X \to X^*$ and $J^*: X^* \to X$ be the normalized duality mappings. We have

(i) $\langle J(x), x \rangle = \|x\| \|J(x)\|_* = \|x\|^2 = \|J(x)\|_*^2$, for any $x \in X$;

(ii) $\langle x^*, J^*(x^*) \rangle = \|J^*(x^*)\| \|x^*\|_* = \|J^*(x^*)\|^2 = \|x^*\|_*^2$, for any $x^* \in X^*$.

The normalized duality mapping $J$ in uniformly convex and uniformly smooth Banach space has many useful properties (see [1, 11, 23]). The following property will be used in this paper.

**Lemma 2.1 in [13].** *Let $X$ be a uniformly convex and uniformly smooth Banach space. For any $x, y \in X$, one has*

$$2\langle J(y), x - y\rangle \leq \|x\|^2 - \|y\|^2 \leq 2\langle J(x), x - y\rangle.$$

Let $C$ be a nonempty closed and convex subset of this uniformly convex and uniformly smooth Banach space $X$. Let $P_C: X \to C$ denote the (standard) metric projection operator. For any $x \in X$, $P_C x \in C$ satisfies

$$\|x - P_C x\| \leq \|x - z\|, \text{ for all } z \in C.$$

The metric projection operator $P_C$ has many useful properties (see [1, 11, 23] for more details).

**Proposition 2.6 in [1]** *Let $X$ be a uniformly convex and uniformly smooth Banach space and $C$ a nonempty closed and convex subset of $X$. Then the metric projection $P_C: X \to C$ satisfies the following properties.*

(i) *The operator $P_C$ is fixed on $C$; that is, $P_C(x) = x$, for any $x \in C$;*

(ii) *$P_C$ has the basic variational properties. For any $x \in X$ and $u \in C$, we have*

$$u = P_C(x) \iff \langle J(x - u), u - z\rangle \geq 0, \text{ for all } z \in C;$$

(iii) *$P_C$ is uniformly continuous on each bounded subset in $X$.*

By the continuity of the single-valued mapping $P_C$, the Mordukhovich derivatives for $P_C$ in uniformly convex and uniformly smooth Banach spaces can be reformulated as follows.

**Definition 1.32 in [17].** Let $X$ be a uniformly convex and uniformly smooth Banach space and $C$ a nonempty closed and convex subset of $X$. The Mordukhovich derivative of the single-valued mapping $P_C$ at $(x, P_C(x))$ is denoted by $\widehat{D}^* P_C(x, P_C(x))$ that is a set-valued mapping $\widehat{D}^* P_C(x, P_C(x)): X^* \rightrightarrows X^*$. It is defined, for any $y^* \in X^*$, by

$$\widehat{D}^* P_C(x, P_C(x))(y^*) := \widehat{D}^* P_C(x)(y^*) = \left\{z^* \in X^*: \limsup_{u \to x} \frac{\langle z^*, u-x\rangle - \langle y^*, P_C(u) - P_C(x)\rangle}{\|u-x\| + \|P_C(u) - P_C(x)\|} \leq 0\right\}.$$

It follows that, if $u \neq P_C(x)$, then

$$\widehat{D}^* P_C(x, u)(y^*) = \emptyset, \text{ for all } y^* \in X^*.$$

By using the Mordukhovich derivatives of $P_C$, the covering constant for $P_C$ at a point is defined (see (2.6) in [2]). For any $\bar{x} \in X$, let $\bar{y} = P_C(\bar{x})$, that is, $(\bar{x}, \bar{y}) \in \text{gph} P_C$. The covering constant for $P_C$ at point $(\bar{x}, \bar{y})$ is defined by

$$\hat{\alpha}(P_C, \bar{x}, \bar{y}) := \sup_{\eta > 0} \inf\{\|z^*\|_*: z^* \in \widehat{D}^* P_C(x, y)(w^*), x \in \mathbb{B}(\bar{x}, \eta), y \in P_C(x) \cap \mathbb{B}(\bar{y}, \eta), \|w^*\|_* = 1\}. \quad (2.1)$$

Since $P_C$ is a single-valued mapping. If we write $y = P_C(x)$, for any $x \in X$, then, (2.1) turns to

$$\hat{\alpha}(P_C, \bar{x}, \bar{y}) = \sup_{\eta>0} \inf\{\|z^*\|_*: z^* \in \widehat{D}^*P_C(x, y)(w^*), x \in \mathbb{B}(\bar{x}, \eta), y \in \mathbb{B}(\bar{y}, \eta), \|w^*\|_* = 1\}. \quad (2.2)$$

## 3. Covering Constants of Metric Projection onto Closed Balls in Uniformly Convex and Uniformly Smooth Banach Spaces

In this section, we consider the metric projection operator onto closed balls in uniformly convex and uniformly smooth Banach spaces. For the properties of the metric projection, the readers could see [1, 11, 23]. Let $(X, \|\cdot\|)$ be a real uniformly convex and uniformly smooth Banach space with topological dual space $(X^*, \|\cdot\|_*)$. Let $\mathbb{B}$ and $\mathbb{B}^*$ denote the unit closed balls in $X$ and $X^*$, respectively. Then, for any $r > 0$, $r\mathbb{B}$ and $r\mathbb{B}^*$ are closed balls with radius $r$ and centered at the origins in $X$ and $X^*$, respectively. Let $\mathbb{S}$ be the unit sphere in $X$. Then, $r\mathbb{S}$ is the sphere in $X$ with radius $r$ and centered at $\theta$. For any $c \in X$ and $r > 0$, let $\mathbb{B}(c, r)$ denote the closed ball in $X$ with radius $r$ and centered at $c$. It follows that $\mathbb{B}(\theta, 1) = \mathbb{B}$ and $\mathbb{B}(\theta, r) = r\mathbb{B}$. The metric projection $P_{r\mathbb{B}}: X \to r\mathbb{B}$ has the following analytic representation (see [1, 11]).

$$P_{r\mathbb{B}}(x) = \begin{cases} x, & \text{for any } x \in r\mathbb{B}, \\ \dfrac{r}{\|x\|} x, & \text{for any } x \notin r\mathbb{B}. \end{cases} \quad (3.1)$$

The following notations play important roles in Fréchet derivatives and Mordukhovich derivatives of $P_{r\mathbb{B}}$. To proceed with the following theorems, we review these notations. For any $x \in r\mathbb{S}$, two subsets $x_r^\uparrow$ and $x_r^\downarrow$ of $X\setminus\{\theta\}$ are defined by

(a) $x_r^\uparrow = \{v \in X\setminus\{\theta\}: \text{there is } \delta > 0 \text{ such that } \|x + tv\| > r, \text{ for all } t \in (0, \delta)\}$;
(b) $x_r^\downarrow = \{v \in X\setminus\{\theta\}: \text{there is } \delta > 0 \text{ such that } \|x + tv\| \leq r, \text{ for all } t \in (0, \delta)\}$.

By using the continuity of the single-valued metric projection $P_{r\mathbb{B}}$ in uniformly convex and uniformly smooth Banach spaces, by Definitions 1.13 and 1.32 in [17], the Mordukhovich derivative of $P_{r\mathbb{B}}$ is formulated. For any $x \in X$, the Mordukhovich derivative of $P_{r\mathbb{B}}$ at the point $(x, P_{r\mathbb{B}}(x))$ is a set-valued mapping $\widehat{D}^*P_{r\mathbb{B}}(x, P_{r\mathbb{B}}(x)): X^* \rightrightarrows X^*$, which is calculated as follows. For any $w^* \in X^*$, we have

$$\widehat{D}^*P_{r\mathbb{B}}(x, P_{r\mathbb{B}}(x))(w^*)$$
$$:= \widehat{D}^*P_{r\mathbb{B}}(x)(w^*)$$
$$= \left\{z^* \in X^*: \limsup_{u \to x} \frac{\langle z^*, u-x \rangle - \langle w^*, P_{r\mathbb{B}}(u) - P_{r\mathbb{B}}(x) \rangle}{\|u-x\| + \|P_{r\mathbb{B}}(u) - P_{r\mathbb{B}}(x)\|} \leq 0\right\}. \quad (3.2)$$

**Proposition 3.1 in [15]**. *Let $r > 0$. For any $x \in X\setminus(r\mathbb{B})$, we have*

$$\theta^* \in \widehat{D}^*P_{r\mathbb{B}}(x)(\lambda J(x)), \text{ for any } \lambda \leq 0.$$

**Theorem 3.2 in [15]**. *Let $X$ be a uniformly convex and uniformly smooth Banach space. For $r > 0$ and for $x \in X$, the Mordukhovich derivative of the metric projection $P_{r\mathbb{B}}: X \to r\mathbb{B}$ at point $x$ has the following properties.*

(i) For every $x \in r\mathbb{B}^o$, we have
$$\widehat{D}^* P_{r\mathbb{B}}(x)(w^*) = \{w^*\}, \text{ for every } w^* \in X^*.$$

(ii) For every $x \in X \setminus r\mathbb{B}$, we have
$$\widehat{D}^* P_{r\mathbb{B}}(x)(w^*) = \left\{ \frac{r}{\|x\|} \left( w^* - \frac{\langle w^*, x \rangle}{\|x\|^2} J(x) \right) \right\}, \text{ for every } w^* \in X^*.$$

In particular, we have

(a) $\widehat{D}^* P_{r\mathbb{B}}(x)(w^*) = \left\{ \frac{r}{\|x\|} w^* \right\}$, if $w^* \perp x$;

(b) $\widehat{D}^* P_{r\mathbb{B}}(x)(\lambda J(x)) = \{\theta^*\}$, for every $\lambda \in \mathbb{R}$.

(iii) If $x \in r\mathbb{S}$, then

(a) $\widehat{D}^* P_{r\mathbb{B}}(x)(\theta^*) = \{\theta^*\}$;
(b) For any $w^* \in X^* \setminus \{\theta^*\}$, we have

Part 1. $\theta^* \in \widehat{D}^* P_{r\mathbb{B}}(x)(w^*)$

$\Rightarrow \langle J(x), -J^*(w^*) \rangle > 0, \langle w^*, x \rangle < 0, \langle J(x), -J^*(w^*) \rangle \langle w^*, x \rangle + r^2 \|w^*\|_*^2 = 0$ and $-J^*(w^*) \notin x_r^!$.

In particular, $\theta^* \notin \widehat{D}^* P_{r\mathbb{B}}(x)(\lambda J(x))$, for any $\lambda > 0$;

Part 2. $w^* = \frac{\langle w^*, x \rangle}{\|x\|^2} x^*$ and $\langle w^*, x \rangle < 0 \Rightarrow \theta^* \in \widehat{D}^* P_{r\mathbb{B}}(x)(w^*)$.

In particular, $\theta^* \in \widehat{D}^* P_{r\mathbb{B}}(x)(\lambda J(x))$, for any $\lambda < 0$.

(c) $\widehat{D}^* P_{r\mathbb{B}}(x)(J(x)) = \emptyset$.

For any $\bar{x} \in X$, let $\bar{y} = P_{r\mathbb{B}}(\bar{x})$. Replace $C$ by $r\mathbb{B}$ in (2.1) and (2.2), then the covering constant for $P_{r\mathbb{B}}$ at $(\bar{x}, \bar{y})$ is defined by

$$\hat{\alpha}(P_{r\mathbb{B}}, \bar{x}, \bar{y}) = \sup_{\eta > 0} \inf\{\|z^*\|_* : z^* \in \widehat{D}^* P_{r\mathbb{B}}(x, y)(w^*), x \in \mathbb{B}(\bar{x}, \eta), y \in \mathbb{B}(\bar{y}, \eta), \|w^*\|_* = 1\}.$$

Now, we prove the main theorem in this section.

**Theorem 3.1**. *Let $X$ be a uniformly convex and uniformly smooth Banach space. Let $r > 0$. For $\bar{x} \in X$ with $\bar{y} = P_{r\mathbb{B}}(\bar{x})$, the covering constant for the metric projection $P_{r\mathbb{B}}$ at $(\bar{x}, \bar{y})$ satisfies*

(a) $\hat{\alpha}(P_{r\mathbb{B}}, \bar{x}, \bar{y}) = 1$, for any $\bar{x} \in r\mathbb{B}^o$;
(b) $\hat{\alpha}(P_{r\mathbb{B}}, \bar{x}, \bar{y}) = 0$, for any $\bar{x} \in X \setminus (r\mathbb{B}^o)$.

*Proof*. Proof of (a). Let $\bar{x} \in r\mathbb{B}^o$ with $\bar{y} = P_{r\mathbb{B}}(\bar{x}) = \bar{x}$. The proof of part (a) is divided into three cases with respect to $\eta > 0$.

Case 1. $0 < \eta < r - \|\bar{x}\|$. In this case, we have

$$\mathbb{B}(\bar{x},\eta) = \mathbb{B}(\bar{y},\eta) \subset r\mathbb{B}^o.$$

This implies that, for any $x \in \mathbb{B}(\bar{x},\eta)$, we have

$$y := P_{r\mathbb{B}}(x) = x \in \mathbb{B}(\bar{x},\eta) = \mathbb{B}(\bar{y},\eta).$$

By part (i) in Theorem 3.2 in [15], for any $x \in \mathbb{B}(\bar{x},\eta)$, we have $y := P_{r\mathbb{B}}(x) = x \in \mathbb{B}(\bar{y},\eta)$. The Mordukhovich derivative of the metric projection $P_{r\mathbb{B}}$ at $(x,y)$ satisfies

$$\widehat{D}^*P_{r\mathbb{B}}(x,y)(w^*) = \{w^*\}, \text{ for every } w^* \in X^*.$$

This implies that, for any fixed $0 < \eta < r - \|\bar{x}\|$, we obtain

$$\inf\{\|z^*\|_*: z^* \in \widehat{D}^*P_{r\mathbb{B}}(x,y)(w^*), x \in \mathbb{B}(\bar{x},\eta), y \in \mathbb{B}(\bar{y},\eta), \|w^*\|_* = 1\}$$

$$= \inf\{\|w^*\|_*: \{w^*\} = \widehat{D}^*P_{r\mathbb{B}}(x,y)(w^*), x \in \mathbb{B}(\bar{x},\eta), y \in \mathbb{B}(\bar{y},\eta), \|w^*\|_* = 1\}$$

$$= 1. \tag{3.3}$$

Case 2. $\eta > r - \|\bar{x}\|$. In this case, $\mathbb{B}(\bar{x},\eta)\setminus r\mathbb{B} \neq \emptyset$. Case 2 is divided into two subcases with respect to $\bar{x} \neq \theta$ or $\bar{x} = \theta$.

Subcase 2.1. $\bar{x} \neq \theta$. In this case, let $u = \frac{\frac{1}{2}(\eta+\|\bar{x}\|+r)}{\|\bar{x}\|}\bar{x}$. By $\eta + \|\bar{x}\| > r$, one has

$$r < \|u\| < \eta + \|\bar{x}\|. \tag{3.4}$$

And

$$\|u - \bar{x}\| = \left\|\frac{\frac{1}{2}(\eta+\|\bar{x}\|+r)}{\|\bar{x}\|}\bar{x} - \bar{x}\right\| = \left|\frac{\frac{1}{2}(\eta+\|\bar{x}\|+r)}{\|\bar{x}\|} - 1\right|\|\bar{x}\| = \left|\frac{1}{2}(\eta + \|\bar{x}\| + r) - \|\bar{x}\|\right|$$

$$= \frac{1}{2}(\eta + \|\bar{x}\| + r) - \|\bar{x}\| = \frac{1}{2}(\eta - \|\bar{x}\| + r) < \eta. \tag{3.5}$$

(3.4) and (3.5) imply $u \in \mathbb{B}(\bar{x},\eta)\setminus r\mathbb{B}$. By (3.1), we have

$$v := P_{r\mathbb{B}}(u) = \frac{r}{\|u\|}u = \frac{r}{\left\|\frac{\frac{1}{2}(\eta+\|\bar{x}\|+r)}{\|\bar{x}\|}\bar{x}\right\|} \cdot \frac{\frac{1}{2}(\eta+\|\bar{x}\|+r)}{\|\bar{x}\|}\bar{x} = \frac{r}{\|\bar{x}\|}\bar{x} \in \partial(r\mathbb{B}).$$

We calculate

$$\|v - \bar{x}\| = \left\|\frac{r}{\|\bar{x}\|}\bar{x} - \bar{x}\right\| = \left|\frac{r}{\|\bar{x}\|} - 1\right|\|\bar{x}\| = r - \|\bar{x}\| < \eta.$$

This implies $v = P_{r\mathbb{B}}(u) \in \mathbb{B}(\bar{x},\eta) = \mathbb{B}(P_{r\mathbb{B}}(\bar{x}),\eta)$.

We can use part (b) of (ii) in Theorem 3.2 in [15] to prove this case. However, we use part (a) of (ii) in Theorem 3.2 in [15] to provide a different approach to prove it. For the given $u \in \mathbb{B}(\bar{x},\eta)\setminus r\mathbb{B}$, by

Khan-Banach Theorem, one can show that there is $k^* \in X^*$ such that

$$\|k^*\|_* = 1 \text{ and } \langle k^*, u \rangle = 0.$$

This means $k^* \perp u$. By (a) of part (ii) in Theorem 3.2 in [15], we have

$$\widehat{D}^* P_{r\mathbb{B}}(u, v)(k^*) = \left\{\frac{r}{\|u\|} k^*\right\}.$$

For the case $\eta > r - \|\bar{x}\|$, this implies that, for such $u \in \mathbb{B}(\bar{x}, \eta), v \in \mathbb{B}(\bar{y}, \eta)$ with $\|u\| > r$, we have

$$\inf\{\|z^*\|_*: z^* \in \widehat{D}^* P_{r\mathbb{B}}(x, y)(w^*), x \in \mathbb{B}(\bar{x}, \eta), y \in \mathbb{B}(\bar{y}, \eta), \|w^*\|_* = 1\}$$

$$\leq \inf\{\|z^*\|_*: z^* \in \widehat{D}^* P_{r\mathbb{B}}(u, v)(w^*), u \in \mathbb{B}(\bar{x}, \eta), v \in \mathbb{B}(\bar{y}, \eta), \|w^*\|_* = 1\}$$

$$\leq \inf\left\{\left\|\frac{r}{\|u\|} k^*\right\|_*: \left\{\frac{r}{\|u\|} k^*\right\} = \widehat{D}^* P_{r\mathbb{B}}(u, v)(k^*), u \in \mathbb{B}(\bar{x}, \eta), v \in \mathbb{B}(\bar{y}, \eta), \|k^*\|_* = 1\right\}$$

$$= \left\|\frac{r}{\|u\|}\right\|$$

$$< 1, \text{ for all } \eta > r - \|\bar{x}\|, \text{ with } \bar{x} \in (r\mathbb{B}^o)\setminus\{\theta\}. \tag{3.8}$$

Subcase 2.2. $\eta > r - \|\bar{x}\|$ and $\bar{x} = \theta$. We use part (b) of (ii) in Theorem 3.2 in [15] to prove a result stronger than (3.8) as follows. In this case, take arbitrarily $u \in X$ with $r < \|u\| < \eta$. Then,

$$v := P_{r\mathbb{B}}(u) = \frac{r}{\|u\|} u \in r\mathbb{S} \subseteq \mathbb{B}(\theta, \eta) = \mathbb{B}(\bar{x}, \eta).$$

By (b) of part (ii) in Theorem 3.2 in [15], we have

$$\widehat{D}^* P_{r\mathbb{B}}(u, v)(\lambda J(u)) = \{\theta^*\}, \text{ for every } \lambda < 0.$$

Let $\lambda = -\frac{1}{\|J(u)\|_*} = -\frac{1}{\|u\|}$. Then $\left\|-\frac{1}{\|J(u)\|_*} J(u)\right\|_* = 1$. For such $u \in \mathbb{B}(\theta, \eta), v \in \mathbb{B}(\theta, \eta)$ with $\eta > \|u\| > r$, we have

$$\inf\{\|z^*\|_*: z^* \in \widehat{D}^* P_{r\mathbb{B}}(x, y)(w^*), x \in \mathbb{B}(\theta, \eta), y \in \mathbb{B}(\theta, \eta), \|w^*\|_* = 1\}$$

$$\leq \inf\{\|z^*\|_*: z^* \in \widehat{D}^* P_{r\mathbb{B}}(u, v)(w^*), u \in \mathbb{B}(\theta, \eta), v \in \mathbb{B}(\theta, \eta), \|w^*\|_* = 1\}$$

$$\leq \inf\left\{\|\theta^*\|_*: \{\theta^*\} = \widehat{D}^* P_{r\mathbb{B}}(u, v)\left(-\frac{1}{\|J(u)\|_*} J(u)\right), u \in \mathbb{B}(\bar{x}, \eta), v \in \mathbb{B}(\bar{y}, \eta), \left\|-\frac{1}{\|J(u)\|_*} J(u)\right\|_* = 1\right\}$$

$$= 0. \tag{3.7}$$

By (3.6) and (3.7), for any $\eta > r - \|\bar{x}\|$, we have

$$\inf\{\|z^*\|_*: z^* \in \widehat{D}^* P_{r\mathbb{B}}(x, y)(w^*), x \in \mathbb{B}(\bar{x}, \eta), y \in \mathbb{B}(\bar{y}, \eta), \|w^*\|_* = 1\} < 1. \tag{3.8}$$

Case 3. $\eta = r - \|\bar{x}\|$. Case 3 is divided into two subcases with respect to $\bar{x} \neq \theta$ or $\bar{x} = \theta$.

Subcase 3.1. We suppose that $\bar{x} \neq \theta$. In this case, since $X$ is a uniformly convex and uniformly smooth Banach space, then, $\mathbb{B}(\bar{x}, \eta) \cap r\mathbb{S} = \left\{\frac{r}{\|\bar{x}\|}\bar{x}\right\}$. This implies

$$\mathbb{B}(\bar{x},\eta)\backslash(r\mathbb{B}^o) = \left\{\frac{r}{\|\bar{x}\|}\bar{x}\right\} \quad \text{and} \quad \mathbb{B}(\bar{x},\eta) = (\mathbb{B}(\bar{x},\eta) \cap (r\mathbb{B}^o)) \cup \left\{\frac{r}{\|\bar{x}\|}\bar{x}\right\}. \tag{3.9}$$

For any $x \in \mathbb{B}(\bar{x}, \eta) \cap (r\mathbb{B}^o)$, similar to the proof of Case 1, one proves

$$\inf\{\|z^*\|_*: z^* \in \widehat{D}^* P_{r\mathbb{B}}(x,y)(w^*), x \in \mathbb{B}(\bar{x},\eta) \cap (r\mathbb{B}^o), y \in \mathbb{B}(\bar{y},\eta), \|w^*\|_* = 1\}$$

$$= \inf\{\|w^*\|_*: \{w^*\} = \widehat{D}^* P_{r\mathbb{B}}(x,y)(w^*), x \in \mathbb{B}(\bar{x},\eta) \cap (r\mathbb{B}^o), y \in \mathbb{B}(\bar{y},\eta), \|w^*\|_* = 1\}$$

$$= 1, \text{ for any } x \in \mathbb{B}(\bar{x},\eta) \cap (r\mathbb{B}^o). \tag{3.10}$$

Hence, there is only one point left for consideration, which is $\frac{r}{\|\bar{x}\|}\bar{x} \in \mathbb{B}(\bar{x},\eta)\backslash(r\mathbb{B}^o)$. Since $\frac{r}{\|\bar{x}\|}\bar{x} \in r\mathbb{S}$, by the special case of (b) in part (iii) of Theorem 3.2 in [15], we have

$$\theta^* \in \widehat{D}^* P_{r\mathbb{B}}\left(\frac{r}{\|\bar{x}\|}\bar{x}, \frac{r}{\|\bar{x}\|}\bar{x}\right)(\lambda J(\frac{r}{\|\bar{x}\|}\bar{x})), \text{ for any } \lambda < 0.$$

Let $\lambda = -\frac{1}{r}$. By the properties of the normalized duality mapping, this implies

$$\theta^* \in \widehat{D}^* P_{r\mathbb{B}}\left(\frac{r}{\|\bar{x}\|}\bar{x}, \frac{r}{\|\bar{x}\|}\bar{x}\right)(-J(\frac{1}{\|\bar{x}\|}\bar{x})).$$

Notice that

$$\left\|-J(\frac{1}{\|\bar{x}\|}\bar{x})\right\|_* = \left\|J(\frac{1}{\|\bar{x}\|}\bar{x})\right\|_* = \left\|\frac{1}{\|\bar{x}\|}\bar{x}\right\| = 1.$$

With respect to the fixed point $\frac{r}{\|\bar{x}\|}\bar{x} \in \mathbb{B}(\bar{x},\eta)\backslash(r\mathbb{B}^o)$, we calculate

$$\inf\left\{\|z^*\|_*: z^* \in \widehat{D}^* P_{r\mathbb{B}}\left(\frac{r}{\|\bar{x}\|}\bar{x}, \frac{r}{\|\bar{x}\|}\bar{x}\right)(w^*), \frac{r}{\|\bar{x}\|}\bar{x} \in \mathbb{B}(\bar{x},\eta), \frac{r}{\|\bar{x}\|}\bar{x} \in \mathbb{B}(\bar{x},\eta), \|w^*\|_* = 1\right\}$$

$$\leq \inf\left\{\|\theta^*\|_*: \theta^* \in \widehat{D}^* P_{r\mathbb{B}}\left(\frac{r}{\|\bar{x}\|}\bar{x}, \frac{r}{\|\bar{x}\|}\bar{x}\right)(-J(\frac{1}{\|\bar{x}\|}\bar{x})), \frac{r}{\|\bar{x}\|}\bar{x} \in \mathbb{B}(\bar{x},\eta), \frac{r}{\|\bar{x}\|}\bar{x} \in \mathbb{B}(\bar{x},\eta), \left\|-J(\frac{1}{\|\bar{x}\|}\bar{x})\right\|_* = 1\right\}$$

$$= 0. \tag{3.11}$$

By (3.10), (3.11) and (3.9), for any $\bar{x} \in (r\mathbb{B}^o)\backslash\{\theta\}$. For $\eta = r - \|\bar{x}\|$, we obtain

$$\inf\{\|z^*\|_*: z^* \in \widehat{D}^* P_{r\mathbb{B}}(x,y)(w^*), x \in \mathbb{B}(\bar{x},\eta), y \in \mathbb{B}(\bar{y},\eta), \|w^*\|_* = 1\} = 0. \tag{3.12}$$

Subcase 3.2. We suppose $\bar{x} = \theta$ with $\eta = r - \|\bar{x}\| = r$. In this case, $\mathbb{B}(\bar{x},\eta) = \mathbb{B}(\theta,r)$, which has the following partition

$$\mathbb{B}(\bar{x},\eta) = (r\mathbb{B}^o) \cup (r\mathbb{S}).$$

Similar to the proof of (3.3), we obtain

$$\inf\{\|z^*\|_*: z^* \in \widehat{D}^*P_{r\mathbb{B}}(x,y)(w^*), x \in r\mathbb{B}^o, y \in r\mathbb{B}^o, \|w^*\|_* = 1\} = 1. \quad (3.13)$$

Similar to the proof of (3.11), we obtain

$$\inf\{\|z^*\|_*: z^* \in \widehat{D}^*P_{r\mathbb{B}}(x,y)(w^*), x \in r\mathbb{S}, y \in r\mathbb{S}, \|w^*\|_* = 1\} = 0. \quad (3.14)$$

By (3.13) and (3.14), for $\bar{x} = \theta$, for $\eta = r$, we have

$$\inf\{\|w^*\|_*: \{w^*\} = \widehat{D}^*P_{r\mathbb{B}}(x,y)(w^*), x \in \mathbb{B}(\bar{x},\eta), y \in \mathbb{B}(\bar{y},\eta), \|w^*\|_* = 1\}$$

$$= \inf\{\|w^*\|_*: \{w^*\} = \widehat{D}^*P_{r\mathbb{B}}(x,y)(w^*), x \in \mathbb{B}(\theta,r), y \in \mathbb{B}(\theta,r), \|w^*\|_* = 1\}$$

$$= 0. \quad (3.15)$$

By (3.15) and (3.12), for any $\bar{x} \in r\mathbb{B}^o$, for $\eta = r - \|\bar{x}\|$, we obtain

$$\inf\{\|z^*\|_*: z^* \in \widehat{D}^*P_{r\mathbb{B}}(x,y)(w^*), x \in \mathbb{B}(\bar{x},\eta), y \in \mathbb{B}(\bar{y},\eta), \|w^*\|_* = 1\} = 0. \quad (3.16)$$

Summarizing cases 1, 2 and 3, by (3.3), (3.8) and (3.16), for any $\bar{x} \in r\mathbb{B}^o$, by (2.2), we have

$$\hat{\alpha}(P_{r\mathbb{B}}, \bar{x}, \bar{y}) = \sup_{\eta > 0} \inf\{\|z^*\|_*: z^* \in \widehat{D}^*P_{r\mathbb{B}}(x,y)(w^*), x \in \mathbb{B}(\bar{x},\eta), y \in \mathbb{B}(\bar{y},\eta), \|w^*\|_* = 1\}$$

$$= \sup_{r - \|\bar{x}\| > \eta > 0} \inf\{\|z^*\|_*: z^* \in \widehat{D}^*P_{r\mathbb{B}}(x,y)(w^*), x \in \mathbb{B}(\bar{x},\eta), y \in \mathbb{B}(\bar{y},\eta), \|w^*\|_* = 1\}$$

$$= 1.$$

This proves part (a):

$$\hat{\alpha}(P_{r\mathbb{B}}, \bar{x}, \bar{y}) = 1, \text{ for any } \bar{x} \in r\mathbb{B}^o.$$

Proof of (b). $\bar{x} \in X \setminus (r\mathbb{B}^o)$. In this case, $\|\bar{x}\| \geq r$ and $\bar{y} = P_{r\mathbb{B}}(\bar{x}) = \frac{r}{\|\bar{x}\|}\bar{x}$. For any $\eta > 0$, we have $\mathbb{B}(\bar{x},\eta) \cap (X \setminus (r\mathbb{B})) \neq \emptyset$. Let $k = \beta \bar{x}$, for some positive number $\beta$ satisfying $1 < \beta < 1 + \frac{\eta}{\|\bar{x}\|}$. We have

$$k \in \mathbb{B}(\bar{x},\eta) \cap (X \setminus (r\mathbb{B}))$$

By (b) in part (ii) of Theorem 3.2 in [15] (or, by Proposition 3.1 in [15]), we have

$$\{\theta^*\} = \widehat{D}^*P_{r\mathbb{B}}(k)(\lambda J(k)), \text{ for any } \lambda < 0.$$

Let $\lambda = -\frac{1}{\|k\|}$. Then $\{\theta^*\} = \widehat{D}^*P_{r\mathbb{B}}(k)(-\frac{1}{\|k\|}J(k))$. Notice that

$$\left\|-\frac{1}{\|w\|}J(k)\right\|_* = \left\|\frac{1}{\|w\|}k\right\| = 1.$$

Since $P_{r\mathbb{B}}(k) = \frac{1}{\|k\|}k = \frac{1}{\|\beta\bar{x}\|}\beta\bar{x} = \frac{1}{\|\bar{x}\|}\bar{x} = P_{r\mathbb{B}}(\bar{x})$, then, we have $P_{r\mathbb{B}}(k) \in \mathbb{B}(P_{r\mathbb{B}}(\bar{x}),\eta)$. By using

the point $k \in \mathbb{B}(\bar{x}, \eta) \cap (X\setminus(r\mathbb{B}))$ and by $\bar{y} = P_{r\mathbb{B}}(\bar{x}) = \frac{r}{\|\bar{x}\|}\bar{x}$, we calculate

$$\inf\{\|z^*\|_*: z^* \in \widehat{D}^* P_{r\mathbb{B}}(x)(w^*), x \in \mathbb{B}(\bar{x}, \eta), y \in \mathbb{B}(\bar{y}, \eta), \|w^*\|_* = 1\}$$

$$\leq \inf\left\{\|\theta^*\|_*: \theta^* \in \widehat{D}^* P_{r\mathbb{B}}(k)\left(-\frac{1}{\|k\|}J(k)\right), k \in \mathbb{B}(\bar{x}, \eta), P_{r\mathbb{B}}(k) \in \mathbb{B}\left(\frac{r}{\|\bar{x}\|}\bar{x}, \eta\right), \left\|-\frac{1}{\|k\|}J(k)\right\|_* = 1\right\}$$

$= 0$, for any $\eta > 0$.

This implies that, for any $\bar{x} \in X\setminus(r\mathbb{B}^o)$, we have

$$\hat{\alpha}(P_{r\mathbb{B}}, \bar{x}, \bar{y}) = \sup_{\eta>0} \inf\{\|z^*\|_*: z^* \in \widehat{D}^* P_{r\mathbb{B}}(x)(w^*), x \in \mathbb{B}(\bar{x}, \eta), P_{r\mathbb{B}}(x) \in \mathbb{B}(\bar{y}, \eta), \|w^*\|_* = 1\} = 0.$$

That is,

$$\hat{\alpha}(P_{r\mathbb{B}}, \bar{x}, \bar{y}) = 0, \text{ for any } \bar{x} \in X\setminus(r\mathbb{B}^o). \qquad \square$$

In the special case, when we consider $r = \infty$ in Theorem 3.1, then, $r\mathbb{B}$ becomes $X$ and $P_X = I_X$ and we obtain the following result.

**Corollary 3.2**. *Let X be a uniformly convex and uniformly smooth Banach space. Let $I_X$ be the identity mapping in X. The covering constant for $I_X$ at $(\bar{x}, \bar{x})$ satisfies*

$$\hat{\alpha}(I_X, \bar{x}, \bar{x}) = 1, \text{ for any } \bar{x} \in X.$$

## 4. Covering Constants for Metric Projection onto Closed and Convex Cylinder in $l_p$

Let $p$, $q$ be positive numbers satisfying $1 < p, q < \infty$ and $\frac{1}{p} + \frac{1}{q} = 1$. $(l_p, \|\cdot\|_p)$ is the standard real uniformly convex and uniformly smooth Banach space of sequences of real numbers with dual space $(l_q, \|\cdot\|_q)$. The origins of both $l_p$ and $l_q$ are exactly same $\theta = \theta^* = (0, 0, \dots)$. In this section, we consider some closed and convex cylinders in $l_p$ defined in [14,15]. The Fréchet differentiability and the Mordukhovich derivatives of the metric projection onto closed and convex cylinders in $l_p$ have been studied in [14,15]. In this section, we use the results from [15] to find the covering constants for the metric projection onto closed and convex cylinders in $l_p$. Here, all notations about closed and convex cylinders in $l_p$ follow the notations used in [14, 15].

Let $J: l_p \to l_q$ be the normalized duality mapping in $l_p$. It has the following representation. For any point $x = (x_1, x_2, \dots) \in l_p$ with $x \neq \theta$, we have

$$J(x) = \left(\frac{|x_1|^{p-1}\text{sign}(x_1)}{\|x\|_p^{p-2}}, \frac{|x_2|^{p-1}\text{sign}(x_2)}{\|x\|_p^{p-2}}, \dots\right) = \left(\frac{|x_1|^{p-2}x_1}{\|x\|_p^{p-2}}, \frac{|x_2|^{p-2}x_2}{\|x\|_p^{p-2}}, \dots\right). \quad (4.1)$$

Similarly, to (4.1), the representation of the normalized duality mapping $J^*: l_q \to l_p$ is given, for any $y^* = (y_1, y_2, \dots) \in l_q$ with $y \neq \theta$, by

$$J^*(y^*) = \left(\frac{|y_1|^{q-1}\text{sign}(y_1)}{\|y\|_q^{q-2}}, \frac{|y_2|^{q-1}\text{sign}(y_2)}{\|y\|_q^{q-2}}, \dots\right) = \left(\frac{|y_1|^{q-2}y_1}{\|y\|_q^{q-2}}, \frac{|y_2|^{q-2}y_2}{\|y\|_q^{q-2}}, \dots\right). \quad (4.2)$$

Let $\mathbb{N}$ denote the set of all positive integers. Let $M$ be a nonempty subset of $\mathbb{N}$ and let $\bar{M} = \mathbb{N}\setminus M$ denote the complement of $M$. We respectively define two subspaces in $l_p$ and $l_p$, with respect to $M$, which induce the definition of closed and convex cylinders in $l_p$ and $l_p$.

$$l_p^M = \{x = (x_1, x_2, \ldots) \in l_p: x_i = 0, \text{ for all } i \in \bar{M}\},$$

$$l_q^M = \{y = (y_1, y_2, \ldots) \in l_q: y_i = 0, \text{ for all } i \in \bar{M}\}.$$

$l_p^M$ and $l_q^M$ are closed subspaces of $l_p$ and $l_q$, respectively. They are the dual spaces of each other. That is, $(l_p^M)^* = l_q^M$. We define a mapping $(\cdot)_M: l_p \to l_p^M$, for $x = (x_1, x_2, \ldots) \in l_p$, by

$$(x_M)_i = \begin{cases} x_i, & \text{for } i \in M, \\ 0, & \text{for } i \notin M, \end{cases} \text{ for } i \in \mathbb{N}.$$

Similarly, we define a mapping $(\cdot)_{\bar{M}}: l_p \to l_p^{\bar{M}}$, for $x = (x_1, x_2, \ldots) \in l_p$ by

$$(x_{\bar{M}})_i = \begin{cases} x_i, & \text{for } i \in \bar{M}, \\ 0, & \text{for } i \notin \bar{M}, \end{cases} \text{ for } i \in \mathbb{N}.$$

Then, $l_p$ and $l_q$ have the following decomposition

$$x = x_M + x_{\bar{M}}, \text{ for any } x \in l_p.$$

and
$$x^* = x_M^* + x_{\bar{M}}^*, \text{ for any } x^* \in l_q. \tag{4.3}$$

**Lemma 4.1 in [14]**. *Let $M$ be a nonempty subset of $\mathbb{N}$. Then $J$ is the normalized duality mapping from $l_p^M$ to $(l_p^M)^* = l_q^M$. That is, $J(x) \in l_q^M$, for any $x \in l_p^M$.*

More precisely speaking, we have

$$J(x)_M = \frac{\|x_M\|_p^{p-2}}{\|x\|_p^{p-2}} J(x_M), \text{ for any } x \in l_p \setminus \{\theta\}.$$

Let $\mathbb{B}_M$ denote the unit closed ball in $l_p^M$. For any $r > 0$, $r\mathbb{B}_M$ is the closed ball with radius $r$ and center origin in $l_p^M$. Let $\mathbb{S}_M$ be the unit sphere in $l_p^M$. Then, $r\mathbb{S}_M$ is the sphere in $l_p^M$ with radius $r$ and center $\theta_M$. We define

$$\mathbb{C}_M = \{x \in l_p: x_M \in \mathbb{B}_M\}.$$

$\mathbb{C}_M$ is called the cylinder in $l_p$ with base $\mathbb{B}_M$. It is a closed and convex subset in $l_p$. For any $r > 0$, $r\mathbb{C}_M$ is the cylinder in $l_p$ with base $r\mathbb{B}_M$, which is a closed and convex subset of $l_p$ satisfying, for any $x \in l_p$,

$$x \in r\mathbb{C}_M \quad \Leftrightarrow \quad \|x_M\|_p \leq r. \tag{4.4}$$

The central axis of $r\mathbb{C}_M$ is the subset $\{x \in l_p: x_M = \theta\}$. The boundary of $r\mathbb{C}_M$ is denoted by $\partial(r\mathbb{C}_M)$ satisfying

$$\partial(r\mathbb{C}_M) = \{x \in l_p \colon \|x_M\|_p = r\}.$$

**Lemma 4.2 in [14].** *For any $r > 0$, the metric projection $P_{r\mathbb{C}_M} \colon l_p \to r\mathbb{C}_M$ satisfies the following formula.*

$$P_{r\mathbb{C}_M}(x) = \begin{cases} x, & \text{for any } x \in r\mathbb{C}_M, \\ \dfrac{r}{\|x_M\|_p} x_M + x_{\bar{M}}, & \text{for any } x \in l_p \setminus r\mathbb{C}_M. \end{cases} \tag{4.5}$$

For any $x \in \partial(r\mathbb{C}_M)$ with $\|x_M\|_p = r$, two subsets $x_r^{\Uparrow}$ and $x_r^{\Downarrow}$ of $l_p$ are defined by

(a) $x_r^{\Uparrow} = \{v \in l_p \colon \text{there is } \delta > 0 \text{ such that } \|(x+tv)_M\|_p > r, \text{ for all } t \in (0, \delta)\}$;
(b) $x_r^{\Downarrow} = \{v \in l_p \colon \text{there is } \delta > 0 \text{ such that } \|(x+tv)_M\|_p \leq r, \text{ for all } t \in (0, \delta)\}$.

For any $r > 0$, the Mordukhovich derivatives of the metric projection $P_{r\mathbb{C}_M}$ onto closed and convex cylinders $r\mathbb{C}_M$ in $l_p$ have been provided in [15].

**Theorem 4.1 in [15].** *For any $r > 0$, Mordukhovich derivatives of the metric projection $P_{r\mathbb{C}_M} \colon l_p \to r\mathbb{C}_M$ have the following representation.*

(i) *For every $x \in (r\mathbb{C}_M)^\circ$, we have*

$$\widehat{D}^* P_{r\mathbb{C}_M}(x)(w^*) = \{w^*\}, \text{ for every } w^* \in l_q.$$

(ii) *For every point $x \in l_p \setminus r\mathbb{C}_M$, we have*

$$\widehat{D}^* P_{r\mathbb{C}_M}(x)(w^*) = \left\{ \frac{r}{\|x_M\|_p} \left( w_M^* - \frac{\langle w_M^*, x_M \rangle}{\|x_M\|_p^2} J(\bar{x}_M) \right) + w_{\bar{M}}^* \right\}$$

$$= \left\{ \frac{r}{\|x_M\|_p} \left( w_M^* - \frac{\langle w_M^*, x_M \rangle \|x\|_p^{p-2}}{\|x_M\|_p^p} (J(x))_M \right) + w_{\bar{M}}^* \right\}, \text{ for all } w^* \in l_q.$$

*In particular,*

$$\widehat{D}^* P_{r\mathbb{C}_M}(x)(\lambda J(x)) = \lambda J(x)_{\bar{M}}, \text{ for every } \lambda \in \mathbb{R},$$

*and* $\quad \widehat{D}^* P_{r\mathbb{C}_M}(x)(\lambda J(x_M)) = \{\theta^*\}, \text{ for every } \lambda \in \mathbb{R};$

(iii) *For any $x \in \partial(r\mathbb{C}_M)$, we have*

(a) $\widehat{D}^* P_{r\mathbb{C}_M}(x)(\theta^*) = \{\theta^*\}$;
(b) *For any $w^* \in l_q \setminus \{\theta^*\}$,*

$$\theta^* \in \widehat{D}^* P_{r\mathbb{C}_M}(x)(w^*) \iff w_{\bar{M}}^* = \theta^*, \ -(J^*(w^*))_M \notin x_r^{\Downarrow} \text{ and } \langle w_M^*, x_M \rangle = -r\|w_M^*\|_q.$$

*In particular, $\theta^* \in \widehat{D}^* P_{r\mathbb{C}_M}(x)(\lambda J(x)_M)$, for every $\lambda \leq 0$;*

(c) $\widehat{D}^* P_{r\mathbb{C}_M}(x)(J(x)) = \emptyset$.

Notice that when $M = \mathbb{N}$, then $\bar{M} = \mathbb{N}\backslash M = \emptyset$, which implies that $l_p^{\mathbb{N}} = l_p$ and $l_q^{\mathbb{N}} = l_q$. In this case, $r\mathbb{C}_{\mathbb{N}}$ becomes $r\mathbb{B}_{\mathbb{N}} = r\mathbb{B}$, which is the closed ball in $l_p$ with radius $r$ and centered at origin $\theta$. The following theorem studies the covering constant for the metric projection $P_{r\mathbb{C}_M}$ in $l_p$, which is considered as an extension of Theorem 3.1 from closed balls to closed and convex cylinders in $l_p$. Therefore, the ideas, the methods and the techniques used in the proof of the following Theorem 4.1 are extensions of that used in the proof of Theorem 3.1 in the previous section. Now, we state and prove the main theorem in this section.

**Theorem 4.1.** *Let $r > 0$. For $\bar{x} \in l_p$ with $\bar{y} = P_{r\mathbb{C}_M}(\bar{x})$, the covering constant for the metric projection $P_{r\mathbb{C}_M}$ at $(\bar{x}, \bar{y})$ satisfies*

(a) $\hat{\alpha}(P_{r\mathbb{C}_M}, \bar{x}, \bar{y}) = 1$, *for any $\bar{x} \in (r\mathbb{C}_M)^\circ$;*
(b) $\hat{\alpha}(P_{r\mathbb{C}_M}, \bar{x}, \bar{y}) = 0$, *for any $\bar{x} \in l_p \backslash (r\mathbb{C}_M)^\circ$.*

*Proof.* Proof of (a). Let $\bar{x} \in (r\mathbb{C}_M)^\circ$ with $\bar{y} = P_{P_{r\mathbb{C}_M}}(\bar{x}) = \bar{x}$ and $\|\bar{x}_M\|_p < r$. The proof of part (a) is divided into three cases with respect to $\eta > 0$.

Case 1. $0 < \eta < r - \|\bar{x}_M\|_p$ Notice that,

$$\|x_M - \bar{x}_M\|_p = \|(x - \bar{x})_M\|_p \leq \|x - \bar{x}\|_p, \text{ for any } x \in l_p.$$

Then, for any $x \in \mathbb{B}(\bar{x}, \eta)$, we have

$$\|x_M\|_p \leq \|x_M - \bar{x}_M\|_p + \|\bar{x}_M\|_p$$

$$= \|(x - \bar{x})_M\|_p + \|\bar{x}_M\|_p$$

$$\leq \|x - \bar{x}\|_p + \|\bar{x}_M\|_p$$

$$\leq \eta + \|\bar{x}_M\|_p$$

$$< r - \|\bar{x}_M\|_p + \|\bar{x}_M\|_p = r.$$

This implies that, for $\bar{x} \in (r\mathbb{C}_M)^\circ$,

$$\mathbb{B}(\bar{x}, \eta) = \mathbb{B}(\bar{y}, \eta) \subset (r\mathbb{C}_M)^\circ, \text{ for } 0 < \eta < r - \|\bar{x}_M\|_p.$$

By part (i) in Theorem 4.1 in [15], for any $x \in \mathbb{B}(\bar{x}, \eta)$, the Mordukhovich derivative of the metric projection $P_{r\mathbb{C}_M}$ at $(x, P_{r\mathbb{C}_M}(x))$ satisfies

$$\widehat{D}^* P_{r\mathbb{C}_M}(x, y)(w^*) = \{w^*\}, \text{ for every } w^* \in X^*.$$

Notice that, for $\bar{x} \in (r\mathbb{C}_M)^\circ$, we have

$$x \in \mathbb{B}(\bar{x}, \eta) \implies y := P_{r\mathbb{C}_M}(x) = x \in \mathbb{B}(\bar{x}, \eta) = \mathbb{B}(\bar{y}, \eta).$$

This implies that, for any $\bar{x} \in (r\mathbb{C}_M)^\circ$ and for any fixed $0 < \eta < r - \|\bar{x}_M\|_p$, we obtain

$$\inf\{\|z^*\|_q : z^* \in \widehat{D}^* P_{r\mathbb{C}_M}(x,y)(w^*), x \in \mathbb{B}(\bar{x},\eta), y \in \mathbb{B}(\bar{y},\eta), \|w^*\|_q = 1\}$$

$$= \inf\{\|w^*\|_q : \{w^*\} = \widehat{D}^* P_{r\mathbb{C}_M}(x,y)(w^*), x \in \mathbb{B}(\bar{x},\eta), y \in \mathbb{B}(\bar{y},\eta), \|w^*\|_q = 1\}$$

$$= 1.$$

Then, for any $\bar{x} \in (r\mathbb{C}_M)^\circ$, we obtain

$$\sup_{0 < \eta < r - \|\bar{x}_M\|_p} \inf\{\|z^*\|_q : z^* \in \widehat{D}^* P_{r\mathbb{C}_M}(x,y)(w^*), x \in \mathbb{B}(\bar{x},\eta), y \in \mathbb{B}(\bar{y},\eta), \|w^*\|_q = 1\} = 1. \quad (4.7)$$

Case 2. $\eta > r - \|\bar{x}_M\|_p$. In this case, $\mathbb{B}(\bar{x},\eta) \setminus r\mathbb{C}_M \neq \emptyset$. Case 2 is divided into two subcases with respect to $\bar{x}_M \neq \theta$ or $\bar{x}_M = \theta$.

Subcase 2.1. We suppose that $\bar{x}_M \neq \theta$ with $\bar{x} = \bar{x}_M + \bar{x}_{\bar{M}}$. In this case, by $\bar{x} \in (r\mathbb{C}_M)^\circ$, we have $0 < \|\bar{x}_M\|_p < r$. Let $u = \frac{\frac{1}{2}(\eta + \|\bar{x}_M\|_p + r)}{\|\bar{x}_M\|_p} \bar{x}_M + \bar{x}_{\bar{M}}$. Then, $u_M = \frac{\frac{1}{2}(\eta + \|\bar{x}_M\|_p + r)}{\|\bar{x}_M\|_p} \bar{x}_M$, which implies

$$\|u_M\|_p = \left\|\frac{\frac{1}{2}(\eta + \|\bar{x}_M\|_p + r)}{\|\bar{x}_M\|_p} \bar{x}_M\right\| = \frac{1}{2}(\eta + \|\bar{x}_M\|_p + r).$$

By $\eta + \|\bar{x}_M\|_p > r$, one has

$$r < \|u_M\|_p < \eta + \|\bar{x}_M\|_p. \quad (4.8)$$

And

$$\|u - \bar{x}\|_p = \left\|\frac{\frac{1}{2}(\eta + \|\bar{x}_M\|_p + r)}{\|\bar{x}_M\|_p} \bar{x}_M + \bar{x}_{\bar{M}} - (\bar{x}_M + \bar{x}_{\bar{M}})\right\|_p$$

$$= \left\|\frac{\frac{1}{2}(\eta + \|\bar{x}_M\|_p + r)}{\|\bar{x}_M\|_p} \bar{x}_M - \bar{x}_M\right\|_p$$

$$= \left|\frac{\frac{1}{2}(\eta + \|\bar{x}_M\|_p + r)}{\|\bar{x}_M\|_p} - 1\right| \|\bar{x}_M\|_p$$

$$= \left|\frac{1}{2}(\eta + \|\bar{x}_M\|_p + r) - \|\bar{x}_M\|_p\right|$$

$$= \frac{1}{2}(\eta + \|\bar{x}_M\|_p + r) - \|\bar{x}_M\|_p$$

$$= \frac{1}{2}(\eta - \|\bar{x}_M\|_p + r) < \eta. \quad (4.9)$$

By (4.8) and (4.9), we get $u \in \mathbb{B}(\bar{x},\eta) \setminus (r\mathbb{C}_M)$. By (4.5), we have

$$v := P_{r\mathbb{C}_M}(u) = \frac{r}{\|u_M\|_p} u_M + u_{\bar{M}}$$

$$= \frac{r}{\left\|\frac{\frac{1}{2}(\eta+\|\bar{x}_M\|_p+r)}{\|\bar{x}_M\|_p}\bar{x}_M\right\|_p} \cdot \frac{\frac{1}{2}(\eta+\|\bar{x}_M\|_p+r)}{\|\bar{x}_M\|_p}\bar{x}_M + \bar{x}_{\bar{M}}$$

$$= \frac{r}{\|\bar{x}_M\|_p}\bar{x}_M + \bar{x}_{\bar{M}} \in \partial(r\mathbb{C}_M).$$

We calculate

$$\|v - \bar{x}\|_p = \left\|\frac{r}{\|\bar{x}_M\|_p}\bar{x}_M - \bar{x}_M\right\|_p$$

$$= \left|\frac{r}{\|\bar{x}_M\|_p} - 1\right|\|\bar{x}_M\|_p$$

$$= r - \|\bar{x}_M\|_p < \eta.$$

This implies $v = P_{r\mathbb{C}_M}(u) \in \mathbb{B}(\bar{x}, \eta) = \mathbb{B}(P_{r\mathbb{C}_M}(\bar{x}), \eta)$.

For the above given $u \in \mathbb{B}(\bar{x},\eta)\backslash r\mathbb{C}_M \subseteq l_p\backslash r\mathbb{C}_M$, by part (ii) in Theorem 4.1 in [15], we have

$$\widehat{D}^*P_{r\mathbb{C}_M}(u)(\lambda J(u_M)) = \{\theta^*\}, \text{ for every } \lambda \in \mathbb{R}. \tag{4.10}$$

By the assumption in subcase 2.1 that $\bar{x}_M \neq \theta$ and $u = \frac{\frac{1}{2}(\eta+\|\bar{x}_M\|_p+r)}{\|\bar{x}_M\|_p}\bar{x}_M + \bar{x}_{\bar{M}}$, we have $J(u)_M \neq \theta$. Take $\lambda = -\frac{1}{\|J(u_M)\|_q}$. By $u \in \mathbb{B}(\bar{x},\eta)$, $v = P_{r\mathbb{C}_M}(u) \in \mathbb{B}(P_{r\mathbb{C}_M}(\bar{x}),\eta)$, and by (4.10), we estimate

$$\inf\{\|z^*\|_q : z^* \in \widehat{D}^*P_{r\mathbb{C}_M}(x,y)(w^*), x \in \mathbb{B}(\bar{x},\eta), y \in \mathbb{B}(\bar{y},\eta), \|w^*\|_q = 1\}$$

$$\leq \inf\{\|z^*\|_q : \{z^*\} = \widehat{D}^*P_{r\mathbb{C}_M}(u,v)(w^*), u \in \mathbb{B}(\bar{x},\eta), v \in \mathbb{B}(\bar{y},\eta), \|w^*\|_q = 1\}$$

$$\leq \inf\left\{\|\theta^*\|_q : \{\theta^*\} = \widehat{D}^*P_{r\mathbb{C}_M}(u,v)\left(-\frac{1}{\|J(u_M)\|_q}J(u_M)\right), u \in \mathbb{B}(\bar{x},\eta), v \in \mathbb{B}(\bar{y},\eta), \left\|-\frac{1}{\|J(u_M)\|_q}J(u_M)\right\|_q = 1\right\}$$

$$= 0, \text{ for any } \eta > r - \|\bar{x}_M\|_p \text{ with } \bar{x}_M \neq \theta. \tag{4.11}$$

Subcase 2.2. We suppose that $\bar{x}_M = \theta$ with $\bar{x} = \bar{x}_{\bar{M}}$. In this case, we have $\eta > r$. We take an arbitrary $u_M \in l_p^M$ with $r < \|u_M\|_p < \eta$. Define $u = u_M + \bar{x}_{\bar{M}}$. This implies $u \in \mathbb{B}(\bar{x},\eta)\backslash r\mathbb{C}_M \subseteq l_p\backslash r\mathbb{C}_M$. By (4.5), we define

$$v := P_{r\mathbb{C}_M}(u) = \frac{r}{\|u_M\|_p}u_M + u_{\bar{M}} = \frac{r}{\|u_M\|_p}u_M + \bar{x}_{\bar{M}} \in \partial(r\mathbb{C}_M).$$

By the assumption that $\bar{x}_M = \theta$, we calculate

$$\|v - \bar{x}\|_p = \left\|\frac{r}{\|u_M\|_p}u_M + \bar{x}_{\bar{M}} - (\bar{x}_M + \bar{x}_{\bar{M}})\right\|_p = r < \eta.$$

This implies $v = P_{r\mathbb{C}_M}(u) \in \mathbb{B}(\bar{x}, \eta) = \mathbb{B}(P_{r\mathbb{C}_M}(\bar{x}), \eta)$. Similar to the proof of (4.11) and by using (4.10), we can show

$$\inf\{\|z^*\|_q : z^* \in \widehat{D}^* P_{r\mathbb{C}_M}(x,y)(w^*), x \in \mathbb{B}(\bar{x}, \eta), y \in \mathbb{B}(\bar{y}, \eta), \|w^*\|_q = 1\}$$

$$\leq \inf\{\|z^*\|_q : \{z^*\} = \widehat{D}^* P_{r\mathbb{C}_M}(u,v)(w^*), u \in \mathbb{B}(\bar{x}, \eta), v \in \mathbb{B}(\bar{y}, \eta), \|w^*\|_q = 1\}$$

$$\leq \inf\left\{\|\theta^*\|_q : \{\theta^*\} = \widehat{D}^* P_{r\mathbb{C}_M}(u,v)\left(-\frac{1}{\|J(u_M)\|_q} J(u_M)\right), u \in \mathbb{B}(\bar{x}, \eta), v \in \mathbb{B}(\bar{y}, \eta), \left\|-\frac{1}{\|J(u_M)\|_q} J(u_M)\right\|_q = 1\right\}$$

$$= 0, \text{ for any } \eta > r - \|\bar{x}_M\|_p \text{ with } \bar{x}_M = \theta. \tag{4.12}$$

By (4.11) and (4.12), we obtain

$$\sup_{\eta > r - \|\bar{x}_M\|_p} \inf\{\|z^*\|_q : z^* \in \widehat{D}^* P_{r\mathbb{C}_M}(x,y)(w^*), x \in \mathbb{B}(\bar{x}, \eta), y \in \mathbb{B}(\bar{y}, \eta), \|w^*\|_q = 1\} = 0. \tag{4.13}$$

Case 3. $\eta = r - \|\bar{x}_M\|_p$. Case 3 is divided into two subcases with respect to $\bar{x}_M \neq \theta$ or $\bar{x}_M = \theta$.

Subcase 3.1. We suppose that $\bar{x}_M \neq \theta$. In this case, since $X$ is a uniformly convex and uniformly smooth Banach space, then, $\mathbb{B}(\bar{x}, \eta) \cap \partial(r\mathbb{C}_M) = \left\{\frac{r}{\|\bar{x}_M\|_p} \bar{x}_M + \bar{x}_{\bar{M}}\right\}$. Let $u := \frac{r}{\|\bar{x}_M\|_p} \bar{x}_M + \bar{x}_{\bar{M}}$. This implies

$$\mathbb{B}(\bar{x}, \eta) \setminus (r\mathbb{C}_M)^o = \{u\} \text{ and } \mathbb{B}(\bar{x}, \eta) = (\mathbb{B}(\bar{x}, \eta) \cap (r\mathbb{C}_M)^o) \cup \{u\}.$$

For any $x \in \mathbb{B}(\bar{x}, \eta) \cap (r\mathbb{C}_M)^o$ with $y = P_{r\mathbb{C}_M}(x) \in \mathbb{B}(\bar{y}, \eta) = \mathbb{B}(\bar{x}, \eta)$. similar to the proof of Case 1, one proves

$$\inf\{\|z^*\|_q : z^* \in \widehat{D}^* P_{r\mathbb{C}_M}(x,y)(w^*), x \in \mathbb{B}(\bar{x}, \eta) \cap (r\mathbb{C}_M)^o, y \in \mathbb{B}(\bar{y}, \eta), \|w^*\|_q = 1\}$$

$$= \inf\{\|w^*\|_q : \{w^*\} = \widehat{D}^* P_{r\mathbb{C}_M}(x,y)(w^*), x \in \mathbb{B}(\bar{x}, \eta) \cap (r\mathbb{C}_M)^o, y \in \mathbb{B}(\bar{y}, \eta), \|w^*\|_q = 1\}$$

$$= 1, \text{ for any } x \in \mathbb{B}(\bar{x}, \eta) \cap (r\mathbb{C}_M)^o. \tag{4.14}$$

Next, we consider the point $u \in \mathbb{B}(\bar{x}, \eta) \setminus (r\mathbb{C}_M)^o$. Since $u \in \partial(r\mathbb{C}_M)$, by the special case of (b) in part (iii) of Theorem 4.1 in [15], we have

$$\theta^* \in \widehat{D}^* P_{r\mathbb{C}_M}(u)(\lambda J(u)_M), \text{ for every } \lambda \leq 0.$$

Let $v := P_{r\mathbb{C}_M}(u)$. By $u \in \partial(r\mathbb{C}_M)$, we have $v = P_{r\mathbb{C}_M}(u) = u$. Let $\lambda = -\frac{1}{\|J(u)_M\|_q}$. By the properties of the normalized duality mapping, this implies

$$\theta^* \in \widehat{D}^* P_{r\mathbb{B}}(u)\left(-\frac{1}{\|J(u)_M\|_q} J(u)_M\right).$$

Notice that

$$\left\|-\frac{1}{\|J(u)_M\|_q}J(u)_M\right\|_q = 1.$$

With respect to the fixed point $u \in \mathbb{B}(\bar{x},\eta)\setminus(r\mathbb{B}^o)$ with $v = P_{r\mathbb{C}_M}(u) = u$, we calculate

$$\inf\{\|z^*\|_q : z^* \in \widehat{D}^*P_{r\mathbb{C}_M}(u,v)(w^*), u \in \mathbb{B}(\bar{x},\eta), v \in \mathbb{B}(\bar{y},\eta), \|w^*\|_q = 1\}$$

$$\leq \inf\left\{\|\theta^*\|_q : \{\theta^*\} = \widehat{D}^*P_{r\mathbb{C}_M}(u,v)\left(-\frac{1}{\|J(u)_M\|_q}J(u)_M\right), x \in \mathbb{B}(\bar{x},\eta), y \in \mathbb{B}(\bar{y},\eta), \left\|-\frac{1}{\|J(u)_M\|_q}J(u)_M\right\|_q = 1\right\}$$

$$= 0. \tag{4.15}$$

By $\mathbb{B}(\bar{x},\eta) = (\mathbb{B}(\bar{x},\eta) \cap (r\mathbb{C}_M)^o) \cup \{u\}$, (4.14) and (4.15), for $\bar{x}_M \neq \theta$, we obtain

$$\inf\{\|z^*\|_q : z^* \in \widehat{D}^*P_{r\mathbb{C}_M}(u,v)(w^*), u \in \mathbb{B}(\bar{x},\eta), v \in \mathbb{B}(\bar{y},\eta), \|w^*\|_q = 1\} = 0, \text{ for } \eta = r - \|\bar{x}_M\|_p. \tag{4.16}$$

Subcase 3.2. We suppose that $\bar{x}_M = \theta$. This implies $\eta = r - \|\bar{x}_M\|_p = r$. In this subcase, we have $\mathbb{B}(\bar{x},\eta) = \mathbb{B}(\theta,r) \subseteq r\mathbb{C}_M$. Hence, $\mathbb{B}(\bar{x},\eta)$ enjoys the following partition.

$$\mathbb{B}(\bar{x},\eta) = (\mathbb{B}(\bar{x},\eta) \cap (r\mathbb{C}_M)^o) \cup (\mathbb{B}(\bar{x},\eta) \cap \partial(r\mathbb{C}_M)).$$

For any $x \in \mathbb{B}(\bar{x},\eta)$, we have $y := P_{r\mathbb{C}_M}(x) = x \in \mathbb{B}(\bar{x},\eta)$. Then, for any $x \in \mathbb{B}(\bar{x},\eta) \cap (r\mathbb{C}_M)^o$, similar to the proof of (4.7) in case 1, we have

$$\inf\{\|z^*\|_q : z^* \in \widehat{D}^*P_{r\mathbb{C}_M}(x,y)(w^*), x \in \mathbb{B}(\bar{x},\eta) \cap (r\mathbb{C}_M)^o, y \in \mathbb{B}(\bar{y},\eta) \cap (r\mathbb{C}_M)^o, \|w^*\|_q = 1\}$$

$$= \inf\{\|w^*\|_q : \{w^*\} = \widehat{D}^*P_{r\mathbb{C}_M}(x,y)(w^*), x \in \mathbb{B}(\bar{x},\eta) \cap (r\mathbb{C}_M)^o, y \in \mathbb{B}(\bar{y},\eta) \cap (r\mathbb{C}_M)^o, \|w^*\|_q = 1\}$$

$$= 1. \tag{4.17}$$

For $x \in \mathbb{B}(\bar{x},\eta) \cap \partial(r\mathbb{C}_M))$, $y := P_{r\mathbb{C}_M}(x) = x \in \mathbb{B}(\bar{y},\eta) \cap \partial(r\mathbb{C}_M))$. Similar to the proof of (4.15), we have

$$\inf\{\|z^*\|_q : z^* \in \widehat{D}^*P_{r\mathbb{C}_M}(x,y)(w^*), x \in \mathbb{B}(\bar{x},\eta) \cap \partial(r\mathbb{C}_M), y \in \mathbb{B}(\bar{y},\eta) \cap \partial(r\mathbb{C}_M), \|w^*\|_q = 1\} = 0. \tag{4.18}$$

By (4.17) and (4.18), for $\bar{x}_M = \theta$, we obtain

$$\inf\{\|z^*\|_q : z^* \in \widehat{D}^*P_{r\mathbb{C}_M}(x,y)(w^*), x \in \mathbb{B}(\bar{x},\eta), y \in \mathbb{B}(\bar{y},\eta), \|w^*\|_q = 1\} = 0. \tag{4.19}$$

Hence, for any $\bar{x} \in (r\mathbb{C}_M)^o$ with $\bar{y} = P_{r\mathbb{C}_M}(\bar{x}) = \bar{x}$, by (4.7) in case 1, (4.13) in case 2 and (4.19) in case 3, we obtain

$$\hat{\alpha}(P_{r\mathbb{C}_M},\bar{x},\bar{y}) = \sup_{\eta>0} \inf\{\|z^*\|_q : z^* \in \widehat{D}^*P_{r\mathbb{C}_M}(x,y)(w^*), x \in \mathbb{B}(\bar{x},\eta), y \in \mathbb{B}(\bar{y},\eta), \|w^*\|_q = 1\}$$

$$= \sup_{0<\eta<r-\|\bar{x}_M\|_p} \inf\{\|z^*\|_q : z^* \in \widehat{D}^*P_{r\mathbb{C}_M}(x,y)(w^*), x \in \mathbb{B}(\bar{x},\eta), y \in \mathbb{B}(\bar{y},\eta), \|w^*\|_q = 1\}$$

$$= 1.$$

This proves part (a) of this theorem:
$$\hat{\alpha}(P_{r\mathbb{C}_M}, \bar{x}, \bar{y}) = 1, \text{ for any } \bar{x} \in r\mathbb{B}^o \text{ with } \bar{y} = P_{r\mathbb{C}_M}(\bar{x}) = \bar{x}.$$

Proof of (b). Let $\bar{x} = \bar{x}_M + \bar{x}_{\bar{M}} \in l_p \setminus (r\mathbb{C}_M)^o$. Then $\|\bar{x}_M\|_p \geq r$. For any given $\eta > 0$, take $u = \beta \bar{x}_M + \bar{x}_{\bar{M}}$, for some $\beta$ with $1 < \beta < 1 + \frac{\eta}{\|\bar{x}_M\|_p}$. This implies

$$\|u - \bar{x}\|_p = \|\beta \bar{x}_M - \bar{x}_M\|_p = (\beta - 1)\|\bar{x}_M\|_p < \eta.$$

Let $v = P_{r\mathbb{C}_M}(u)$. Then

$$v = P_{r\mathbb{C}_M}(u) = \frac{r}{\|\beta \bar{x}_M\|_p} \beta \bar{x}_M + \bar{x}_{\bar{M}} = \frac{r}{\|\bar{x}_M\|_p} \bar{x}_M + \bar{x}_{\bar{M}} = P_{r\mathbb{C}_M}(\bar{x})$$

This implies that, for $u \in \mathbb{B}(\bar{x}, \eta)$, we have $v = P_{r\mathbb{C}_M}(u) \in \mathbb{B}(P_{r\mathbb{C}_M}(\bar{x}), \eta) = \mathbb{B}(\bar{y}, \eta)$. For the point $u \in l_p \setminus r\mathbb{C}_M$, by part (ii) in Theorem 4.1 in [15], we have

$$\widehat{D}^* P_{r\mathbb{C}_M}(u)(\lambda J(u_M)) = \{\theta^*\}, \text{ for every } \lambda < 0.$$

Take $\lambda = -\frac{1}{\|J(u_M)\|_q}$. We have $\left\| -\frac{1}{\|J(u_M)\|_q} J(u_M) \right\|_q = 1$. For the arbitrarily given $\eta > 0$ with the points $u \in \mathbb{B}(\bar{x}, \eta)$ and $v = P_{r\mathbb{C}_M}(u) \in \mathbb{B}(P_{r\mathbb{C}_M}(\bar{x}), \eta)$, we have

$$\inf\{\|z^*\|_q : z^* \in \widehat{D}^* P_{r\mathbb{C}_M}(x, y)(w^*), x \in \mathbb{B}(\bar{x}, \eta), y \in \mathbb{B}(\bar{y}, \eta), \|w^*\|_q = 1\}$$

$$\leq \inf\{\|z^*\|_q : z^* \in \widehat{D}^* P_{r\mathbb{C}_M}(u, v)(w^*), u \in \mathbb{B}(\bar{x}, \eta), v \in \mathbb{B}(\bar{y}, \eta), \|w^*\|_q = 1\}$$

$$\leq \inf\left\{\|\theta^*\|_q : \{\theta^*\} = \widehat{D}^* P_{r\mathbb{C}_M}(u, v)\left(-\frac{1}{\|J(u_M)\|_q} J(u_M)\right), u \in \mathbb{B}(\bar{x}, \eta), v \in \mathbb{B}(\bar{y}, \eta), \left\|-\frac{1}{\|J(u_M)\|_q} J(u_M)\right\|_q = 1\right\}$$

$$= 0.$$

This implies

$$\hat{\alpha}(P_{r\mathbb{C}_M}, \bar{x}, \bar{y})$$

$$= \sup_{\eta > 0} \inf\{\|z^*\|_q : z^* \in \widehat{D}^* P_{r\mathbb{C}_M}(x, y)(w^*), x \in \mathbb{B}(\bar{x}, \eta), y \in \mathbb{B}(\bar{y}, \eta), \|w^*\|_q = 1\}$$

$$\leq \sup_{\eta > 0} \inf\left\{\|\theta^*\|_q : \{\theta^*\} = \widehat{D}^* P_{r\mathbb{C}_M}(u, v)\left(-\frac{1}{\|J(u_M)\|_q} J(u_M)\right), u \in \mathbb{B}(\bar{x}, \eta), v \in \mathbb{B}(\bar{y}, \eta), \left\|-\frac{1}{\|J(u_M)\|_q} J(u_M)\right\|_q = 1\right\}$$

$$= 0, \text{ for any } \bar{x} \in l_p \setminus (r\mathbb{C}_M)^o. \qquad \square$$

In Theorem 4.1, let $M = \mathbb{N}$. It reduces $r\mathbb{C}_M = r\mathbb{B}$, which is a closed ball in $l_p$. We obtain the following corollary immediately, which is a special case of Theorem 3.1 with $X = l_p$.

**Corollary 4.2.** *Let $r > 0$. For $\bar{x} \in l_p$ with $\bar{y} = P_{r\mathbb{B}}(\bar{x})$, we have*

(a) $\hat{\alpha}(P_{r\mathbb{B}}, \bar{x}, \bar{y}) = 1$, for any $\bar{x} \in r\mathbb{B}^\circ$;

(b) $\hat{\alpha}(P_{r\mathbb{B}}, \bar{x}, \bar{y}) = 0$, for any $\bar{x} \in l_p \backslash (r\mathbb{B}^\circ)$.

## 5. Covering Constants for Metric Projection onto the Positive Cone in $L_p(S)$

In theory of function spaces, one of the most important classes of function spaces are the real uniformly convex and uniformly smooth Banach spaces $(L_p(S), \|\cdot\|_p)$, in which $(S, \mathcal{A}, \mu)$ is a positive and complete measure space and $p > 1$. The dual space of $(L_p(S), \|\cdot\|_p)$ is the real uniformly convex and uniformly smooth Banach space $(L_q(S), \|\cdot\|_q)$, in which $p$ and $q$ satisfy $1 < p, q < \infty$ and $\frac{1}{p} + \frac{1}{q} = 1$. $L_p(S)$ and $L_q(S)$ share the same origin $\theta = \theta^*$. For the clarity of the distinction between $L_p(S)$ and its dual space $L_q(S)$, we use English letters $f, g, h, \ldots$ for the elements in $L_p(S)$, and we use Greek letters $\varphi, \psi, \xi, \ldots$ for the elements in the dual space $L_q(S)$.

We define the positive cone $K_p$ of $L_p(S)$ and the positive cone $K_q$ of $L_q(S)$ as follows:

$$K_p = \{f \in L_p(S): f(s) \geq 0, \text{ for } \mu\text{-almost all } s \in S\}.$$

$$K_q = \{\varphi \in L_q(S): \varphi(s) \geq 0, \text{ for } \mu\text{-almost all } s \in S\}.$$

$K_p$ and $K_q$ are pointed closed and convex cones. Let $\leqslant_p$ and $\leqslant_q$ be the two partial orders on $L_p(S)$ and $L_q(S)$ induced by the pointed closed and convex cones $K_p$ and $K_q$, respectively. For any $f, g \in L_p(S)$ with $f \leqslant_p g$, and for any $\xi, \psi \in L_q(S)$ with $\xi \leqslant_q \psi$, we define the ordered intervals in $L_p(S)$ and $L_q(S)$ by

$$[f, g]_{\leqslant_p} = \{h \in L_p(S): f \leqslant_p h \leqslant_p g\},$$

$$[\xi, \psi]_{\leqslant_q} = \{\varphi \in L_q(S): \xi \leqslant_q \varphi \leqslant_q \psi\}.$$

The normalized duality mapping $J: L_p(S) \to L_q(S)$ has the following representation, for any given $f \in L_p(S)$ with $f \neq \theta$,

$$(Jf)(s) = \frac{|f(s)|^{p-1} \text{sign}(f(s))}{\|f\|_p^{p-2}} = \frac{|f(s)|^{p-2} f(s)}{\|f\|_p^{p-2}}, \text{ for all } s \in S.$$

In [14, 15], the Fréchet differentiability and Mordukhovich derivatives of the metric projection $P_{K_p}: L_p(S) \to K_p$ onto closed and convex positive cone in $L_p(S)$ have been studied. In this section, we use the results from [15] to find the covering constant for the metric projection $P_{K_p}$. Here, we first review some properties of $P_{K_p}$ obtained in [14, 15].

**Lemma 5.1 in [14]**. *Both $K_p$ and $K_q$ have empty interior.*

**Lemma 5.2 in [15]**. *The metric projection $P_{K_p}: L_p(S) \to K_p$ has the following representation.*

$$(P_{K_p} f)(s) = \begin{cases} f(s), & \text{if } f(s) > 0, \\ 0, & \text{if } f(s) \leq 0, \end{cases} \text{ for any } f \in L_p(S).$$

**Lemma 5.3 in [15]**. *The metric projection* $P_{K_p}: L_p(S) \to K_p$ *has the following properties*

(a) $P_{K_p}(f) = f$, *for any* $f \in K_p$;
(b) $P_{K_p}(f) = \theta$, *for any* $f \in -K_p$;
(c) $P_{K_p}(f + g) = f + g$, *for any* $f, g \in K_p$;
(d) *For any* $f \in L_p(S)$,

$$P_{K_p}(\lambda f) = \lambda P_{K_p}(f), \text{ for any } \lambda \geq 0.$$

We need the following notations. For any $f \in L_p(S)$, for all $s \in S$, we write

$$f^+(s) = \begin{cases} f(s), & \text{for } f(s) > 0, \\ 0, & \text{for } f(s) \leq 0. \end{cases}$$

and

$$f^-(s) = \begin{cases} f(s), & \text{for } f(s) < 0, \\ 0, & \text{for } f(s) \geq 0. \end{cases}$$

**Lemma 5.1 in [15]**. *For any* $f \in L_p(S) \setminus \{\theta\}$, $f^+$ *and* $f^-$ *have the following properties.*

(a) $f = f^+ + f^-$;
(b) $P_{K_p}(f) = f^+$;
(c) $f = P_{K_p}(f) + f^-$;
(d) $(J(f))^+ = \dfrac{\|f^+\|_p^{p-2}}{\|f\|_p^{p-2}} J(f^+)$;
(e) $(J(f))^- = \dfrac{\|f^-\|_p^{p-2}}{\|f\|_p^{p-2}} J(f^-)$;
(f) $\langle J(f), f^+ \rangle = \langle (J(f))^+, f^+ \rangle = \dfrac{\|f^+\|_p^{p-2}}{\|f\|_p^{p-2}} \langle J(f^+), f^+ \rangle = \dfrac{\|f^+\|_p^p}{\|f\|_p^{p-2}}$;
(g) $\langle J(f), f^- \rangle = \langle (J(f))^-, f^- \rangle = \dfrac{\|f^-\|_p^{p-2}}{\|f\|_p^{p-2}} \langle J(f^-), f^- \rangle = \dfrac{\|f^-\|_p^p}{\|f\|_p^{p-2}}$.

**Theorem 5.2 in [15]**. *The metric projection* $P_{K_p}: L_p(S) \to K_p$ *has the following Mordukhovich derivatives. For any* $f \in L_p(S)$,

(i) $\widehat{D}^* P_{K_p}(f)(\theta^*) = \{\theta^*\}$;
(ii) *For any* $\varphi \in L_q(S)$,

$$\theta^* \in \widehat{D}^* P_{K_p}(f)(\varphi)$$

$$\Leftrightarrow \mu(\{s \in S: \varphi(s) \neq 0 \text{ and } f(s) > 0\} \cup \{s \in S: \varphi(s) < 0 \text{ and } f(s) \leq 0\}) = 0.$$

*This implies that,*

(b₁) *For any* $f \in -K_p$, *we have*

$$\theta^* \in \widehat{D}^* P_{K_p}(f)(\varphi), \text{ for any } \varphi \in K_q;$$

(b₂) *For any $f \in K_p \setminus \{\theta\}$,*

$$\theta^* \notin \widehat{D}^* P_{K_p}(f)(J(f));$$

(iii) $J(f) \in \widehat{D}^* P_{K_p}(f)(J(f))$, *for any $f \in K_p$;*

(iv) *For any given $\psi \in K_q$, we have*

$$\widehat{D}^* P_{K_p}(\theta)(\psi) = [\theta^*, \psi]_{\leqslant q}.$$

It is well-known that the metric projection operator onto nonempty closed and convex subsets in Hilbert spaces is nonexpansive. However, in uniformly convex and uniformly smooth Banach spaces, the metric projection operator does not enjoy the nonexpansiveness property, in general. But, in $L_p(S)$, we will prove that the metric projection $P_{K_p}: L_p(S) \to K_p$ is indeed nonexpansiveness. This is a very important property, which should be very useful in nonlinear analysis related with the function space $L_p(S)$, such as, approximation theory, fixed point theory, optimization theory, variational analysis, and so forth. In the following proposition, we prove the nonexpansiveness of $P_{K_p}$.

**Proposition 5.1.** *The metric projection $P_{K_p}: L_p(S) \to K_p$ is nonexpansive. That is,*

$$\left\| P_{K_p}(f) - P_{K_p}(g) \right\|_p \leq \|f - g\|_p, \text{ for any } f, g \in L_p(S).$$

Proof. For any $f, g \in L_p(S)$, by part (b) in Lemma 5.1 in [15], we estimate

$$\left\| P_{K_p}(f) - P_{K_p}(g) \right\|_p^p$$

$$= \|f^+ - g^+\|_p^p$$

$$= \int_S |f^+(s) - g^+(s)|^p \mu(ds)$$

$$= \int_{f(s)>0, g(s)>0} |f^+(s) - g^+(s)|^p \mu(ds) + \int_{f(s)>0, g(s)\leq 0} |f^+(s) - g^+(s)|^p \mu(ds)$$

$$+ \int_{f(s)\leq 0, g(s)>0} |f^+(s) - g^+(s)|^p \mu(ds) + \int_{f(s)\leq 0, g(s)\leq 0} |f^+(s) - g^+(s)|^p \mu(ds)$$

$$= \int_{f(s)>0, g(s)>0} |f(s) - g(s)|^p \mu(ds) + \int_{f(s)>0, g(s)\leq 0} |f(s) - 0|^p \mu(ds)$$

$$+ \int_{f(s)\leq 0, g(s)>0} |0 - g(s)|^p \mu(ds) + \int_{f(s)\leq 0, g(s)\leq 0} |0 - 0|^p \mu(ds)$$

$$\leq \int_{f(s)>0, g(s)>0} |f(s) - g(s)|^p \mu(ds) + \int_{f(s)>0, g(s)\leq 0} |f(s) - g(s)|^p \mu(ds)$$

$$+ \int_{f(s)\leq 0, g(s)>0} |f(s) - g(s)|^p \mu(ds) + \int_{f(s)\leq 0, g(s)\leq 0} |f(s) - g(s)|^p \mu(ds)$$

$$= \int_S |f(s) - g(s)|^p \mu(ds)$$

$$= \|f - g\|_p^p. \qquad \square$$

Now, by using the nonexpansiveness of the metric projection $P_{K_p}$, we prove the main theorem in this section.

**Theorem 5.2.** *The covering constant for the metric projection $P_{K_p}: L_p(S) \to K_p$ satisfies*

$$\hat{\alpha}\left(P_{K_p}, \bar{f}, \bar{g}\right) = 0, \text{ for any } \bar{f} \in L_p(S) \text{ with } \bar{g} = P_{K_p}(\bar{f}) \in K_p.$$

Proof. Let $\bar{f} \in L_p(S)$ with $\bar{g} := P_{K_p}(\bar{f}) \in K_p$. Since $K_p \cup (-K_p)$ has empty interior, for any $\eta > 0$, there is $f \in \mathbb{B}(\bar{f}, \eta) \cap (L_p(S) \backslash (K_p \cup (-K_p)))$ satisfying

$$\mu\{s \in S: f(s) > 0\} > 0 \quad \text{and} \quad \mu\{s \in S: f(s) < 0\} > 0.$$

Take $E \subseteq \{s \in S: f(s) < 0\}$ such that $0 < \mu(E) < \infty$. Define $\varphi \in L_q(S)$, for every $s \in S$, by

$$\varphi(s) = \begin{cases} \dfrac{1}{\mu(E)^{\frac{1}{q}}}, & \text{if } s \in E, \\ 0, & \text{if } s \notin E. \end{cases}$$

Then, the above $f \in L_p(S)$ and $\varphi \in L_q(S)$ satisfy the following equation

$$\mu(\{s \in S: \varphi(s) \neq 0 \text{ and } f(s) > 0\} \cup \{s \in S: \varphi(s) < 0 \text{ and } f(s) \leq 0\}) = 0.$$

By part (ii) in Theorem 5.2 in [15], we have

$$\theta^* \in \widehat{D}^* P_{K_p}(f)(\varphi).$$

Since $f \in \mathbb{B}(\bar{f}, \eta) \cap (L_p(S) \backslash (K_p \cup (-K_q)))$, by Proposition 5.1, we have

$$\left\| P_{K_p}(f) - P_{K_p}(\bar{f}) \right\|_p \leq \|f - \bar{f}\|_p < \eta.$$

We obtain that

$$f \in \mathbb{B}(\bar{f}, \eta) \quad \Rightarrow \quad P_{K_p}(f) \in \mathbb{B}\left(P_{K_p}(\bar{f}), \eta\right).$$

For the arbitrarily given $\eta > 0$, by $f \in \mathbb{B}(\bar{f}, \eta) \cap (L_p(S) \backslash (K_p \cup (-K_q)))$ and $\varphi \in L_q(S)$ as given above, we obtain

$$\inf\left\{\|\psi\|_q: \psi \in \widehat{D}^* P_{K_p}\left(h, P_{K_p}(h)\right)(\omega), h \in \mathbb{B}(\bar{f}, \eta), P_{K_p}(h) \in \mathbb{B}\left(P_{K_p}(\bar{f}), \eta\right), \|\omega\|_q = 1\right\}$$

$$\leq \inf\left\{\|\psi\|_q: \psi \in \widehat{D}^* P_{K_p}\left(f, P_{K_p}(f)\right)(\omega), f \in \mathbb{B}(\bar{f}, \eta), P_{K_p}(f) \in \mathbb{B}\left(P_{K_p}(\bar{f}), \eta\right), \|\omega\|_q = 1\right\}$$

$$\leq \inf\left\{\|\psi\|_q : \psi \in \widehat{D}^*P_{K_p}\left(f, P_{K_p}(f)\right)(\varphi), f \in \mathbb{B}(\bar{f},\eta), P_{K_p}(f) \in \mathbb{B}\left(P_{K_p}(\bar{f}),\eta\right), \|\varphi\|_q = 1\right\}$$

$$\leq \inf\left\{\|\theta^*\|_q : \theta^* \in \widehat{D}^*P_{K_p}\left(f, P_{K_p}(f)\right)(\varphi), f \in \mathbb{B}(\bar{f},\eta), P_{K_p}(f) \in \mathbb{B}\left(P_{K_p}(\bar{f}),\eta\right), \|\varphi\|_q = 1\right\}$$

$$= 0.$$

This implies

$$\hat{\alpha}(P_{K_p}, \bar{f}, P_{K_p}(\bar{f}))$$

$$= \sup_{\eta>0} \inf\left\{\|\psi\|_q : \psi \in \widehat{D}^*P_{K_p}\left(h, P_{K_p}(h)\right)(\omega), h \in \mathbb{B}(\bar{f},\eta), P_{K_p}(h) \in \mathbb{B}\left(P_{K_p}(\bar{f}),\eta\right), \|\omega\|_q = 1\right\}$$

$$= 0, \text{ for any } \bar{f} \in L_p(S). \qquad \square$$

### 6. Stochastic Fixed-Point Problems

The result of Theorem 3.1 in [2] is very strong. It proves an existence theorem for coincidence points of parameterized set-valued mappings (multifunctions), which is immediately applied to the existence of solutions of parameterized generalized equations, implicit functions, fixed-point theorems, optimal value functions in parametric optimization. In this section, we first find the covering constants of some linear and continuous operators in Banach spaces, which can be considered as special cases of the metric projection operator. Then, we apply Theorem 3.1 in [2] to solving some stochastic fixed-point problems associated with both set-valued mappings and single-valued mappings. Three examples for solving stochastic fixed-point problems are proved in the end of this section.

Before we proceed to investigate the solvability of stochastic fixed-point problems, we need to prove a proposition about the Fréchet differentiability of linear and continuous operators in Banach spaces. For the sake of convenience, we review the definition of Fréchet differentiability of single-valued mappings in Banach spaces.

**Definition 1.13 in [17]**. Let $X$ be a Banach space and let $f: X \to X$ be a single-valued mapping. Let $x \in X$. If there is a linear and continuous mapping $\nabla f(x): X \to X$ such that

$$\lim_{u \to x} \frac{f(u) - f(x) - \nabla f(x)(u-x)}{\|u-x\|} = \theta.$$

then $f$ is said to be Fréchet differentiable at $x$; $\nabla f(x)$ is called the Fréchet derivative of $f$ at $x$. The following theorem provides the connection between Fréchet derivatives and Mordukhovich derivatives.

**Theorem 1.38 in [17]**. *Let $X$ be a Banach space with dual space $X^*$ and let $f: X \to X$ be a single-valued mapping. Suppose that $f$ is Fréchet differentiable at $x \in X$ with $y = f(x)$. Then, the Mordukhovich derivative of $f$ at $x$ satisfies the following equation*

$$\widehat{D}^*f(x,y)(y^*) = \{(\nabla f(x))^*(y^*)\}, \text{ for all } y^* \in X^*.$$

The following proposition is an extension of Corollary 3.2 in section 3 of this paper. The results of the following proposition will be used to prove the main theorem of this section. Note that, in the proposition, theorem and corollaries in this section, the underlying Banach space is not required to be uniformly convex and uniformly smooth.

**Proposition 6.1**. *Let $(X, \|\cdot\|)$ be a real Banach space with dual space $(X^*, \|\cdot\|_*)$. Let $I_X$ be the identity mapping in X. For any real number $\lambda$, the linear and continuous mapping $\lambda I_X: X \to X$ satisfies*

(i) $\lambda I_X$ *is Fréchet differentiable at every point in X such that*

$$\nabla(\lambda I_X)(x) = \lambda I_X, \text{ for any } x \in X;$$

(ii) *The Mordukhovich derivative of $\lambda I_X$ is*

$$\widehat{D}^*(\lambda I_X)(x, \lambda x) = \lambda I_{X^*}, \text{ for any } x \in X.$$

*That is, for any $x \in X$, $\widehat{D}^*(\lambda I_X)$ has the following representation*

$$\widehat{D}^*(\lambda I_X)(x, \lambda x)(y^*) = \{\lambda y^*\}, \text{ for all } y^* \in X^*.$$

(iii) *In addition, if $|\lambda| \leq 1$, then the covering constant for $\lambda I_X$ is constant in X with*

$$\widehat{\alpha}(\lambda I_X, x, \lambda x) = |\lambda|, \text{ for any } x \in X.$$

Proof. Proof of (i). We have

$$\lim_{u \to x} \frac{\lambda I_X(u) - \lambda I_X(x) - \lambda I_X(u-x)}{\|u-x\|} = \theta, \text{ for any given } x \in X.$$

By Definition 1.13 in [17], this proves (i). By Theorem 1.38 in [17], part (ii) follows from (i) immediately. Now, by (ii), we prove (iii). For any $x \in X$, if $|\lambda| \leq 1$, then, for any $\eta > 0$, we have

$$u \in \mathbb{B}(x, \eta) \implies \lambda u \in \mathbb{B}(\lambda x, \eta), \text{ for any } u \in X.$$

We calculate the covering constant for $\lambda I_X$ at point $x$.

$$\widehat{\alpha}(\lambda I_X, x, \lambda x) = \sup_{\eta > 0} \inf\{\|\psi\|_*: \psi \in \widehat{D}^*(\lambda I_X)(u, \lambda u)(y^*), u \in \mathbb{B}(x, \eta), \lambda u \in \mathbb{B}(\lambda x, \eta), \|y^*\|_* = 1\}$$

$$= \sup_{\eta > 0} \inf\{\|\lambda y^*\|_*: \{\lambda y^*\} = \widehat{D}^*(\lambda I_X)(u, \lambda u)(y^*), u \in \mathbb{B}(x, \eta), \lambda u \in \mathbb{B}(\lambda x, \eta), \|y^*\|_* = 1\}$$

$$= |\lambda|, \text{ for any } x \in X. \qquad \square$$

Let $X$ be a Banach space and $U$ and $V$ nonempty subsets in $X$ and let $\mathbb{B}_X$ be the unit ball in $X$. Let $G: X \rightrightarrows X$ be a set-valued mapping. $G$ is said to be Lipschitz-like on $U$ relative to $V$ with some modulus $l \geq 0$ if we have the inclusion (see (2.3) in [2])

$$G(x) \cap V \subset G(u) + l\|x - u\|\mathbb{B}_X, \text{ for all } x, u \in U.$$

In particular, let $V = X$, then, $G$ is said to be Lipschitz-like on $U$ with some modulus $l \geq 0$ if we have the inclusion (see (2.3) in [2])

$$G(x) \subset G(u) + l\|x - u\|\mathbb{B}_X, \text{ for all } x, u \in U.$$

Let $(S, \tau, \mu)$ be a topological probability space, in which $S$ is the sample space. Meanwhile, $(S, \tau)$ is a topological space with the topology $\tau$ on $S$ that coincides with the $\sigma$-field of all events in $S$, and $\mu$ is the probability measure in $S$ defined on $\tau$. Let $G: X \times S \rightrightarrows X$ be a set-valued mapping. Let $(x, s) \in X \times S$. If it satisfies

$$x \in G(x, s),$$

then, $x$ is called a stochastic fixed-point of $G$ with respect to the possible outcome $s \in S$. The *stochastic fixed-point problem* associated with the set-valued mapping $G$ with respect to the probability space $S$ is to find a fixed point of $G$ with respect to a possible outcome in $S$.

Furthermore, let $\sigma: S \to X$ be a single-valued mapping. If $\sigma$ is $\tau$-measurable, then $\sigma$ could be considered as a generalized stochastic variable (or a generalized random variable) on $S$ with values in $X$. Let $\bar{s} \in S$. If there is a neighborhood (an event) $W \subset S$ of $\bar{s}$ such that

$$\sigma(s) \in G(\sigma(s), s), \text{ for any } s \in W,$$

then $\sigma$ is called a *stochastic fixed-point variable* on $W$ around the possible outcome $\bar{s}$. Here, although $\sigma$ is unnecessary to be $\tau$-measurable, but we still use the term stochastic fixed-point variable for $\sigma$. However, in the examples in this section, all stochastic fixed-point variables are indeed $\tau$-measurable.

The underlying Banach space in Theorem 3.1 in [2] is Asplund. We review the concept of Asplund Banach space here. A Banach space $Z$ is Asplund if every convex continuous function defined on an open convex set $O$ in $Z$ is Fréchet differentiable on a dense subset of $O$. The class of Asplund Banach spaces is very large including all reflexive Banach spaces. This implies that every uniformly convex and uniformly smooth Banach space is Asplund; and therefore, the underlying Banach space in Theorem 6.2 is more general than the class of uniformly convex and uniformly smooth Banach spaces. As an application of Theorem 3.1 in [2], we prove the following theorem.

**Theorem 6.2.** *Let X be an Asplund Banach space and let $(S, \tau, \mu)$ be a topological probability space. Let $G: X \times S \rightrightarrows X$ be a set-valued mapping. Suppose that $G$ satisfies the following conditions with respect to a point $(\bar{x}, \bar{s}) \in X \times S$*

(A2)′ *There are neighborhoods U, V of $\bar{x}$ in X, and O of $\bar{s} \in S$ as well as a number l with $0 < l < 1$ such that the multifunction $G(\cdot, s)$ is Lipschitz-like on U relative to V for each $s \in O$ with the uniform modulus l, while the multifunction $s \to G(\bar{x}, s)$ is lower/inner semicontinuous at $\bar{s}$.*

*Then, for any $\lambda, \alpha$ with $1 \geq \lambda > \alpha > l$, there exist a neighborhood $W_{\lambda\alpha} \subset S$ of $\bar{s}$ and a single-valued mapping $\sigma_{\lambda\alpha}: W_{\lambda\alpha} \to X$ such that*

$$\lambda\sigma_{\lambda\alpha}(s) \in G(\sigma_{\lambda\alpha}(s), s), \text{ for any } s \in W_{\lambda\alpha}, \tag{6.1}$$

and

$$\|\sigma_{\lambda\alpha}(s) - \bar{x}\| \leq \frac{\text{dist}(\lambda\bar{x}, G(\bar{x}, s))}{\alpha - l}, \quad \text{for any } s \in W_{\lambda\alpha}. \tag{6.2}$$

*Proof.* In Theorem 3.1 in [2], let $Y = X$. For any given $\lambda, \alpha$ with $1 \geq \lambda > \alpha > l$, we substitute the set-valued mapping $F$ in Theorem 3.1 in [2] by the single-valued linear and continuous mapping $\lambda I_X$. By Proposition 6.1, the covering constant for $\lambda I_X$ at point $(\bar{x}, \lambda\bar{x})$ satisfies

$$\hat{\alpha}(\lambda I_X, \bar{x}, \lambda\bar{x}) = \lambda, \text{ for any } \bar{x} \in X. \tag{6.3}$$

Since the graph of $\lambda I_X$ is the line $\{(x, \lambda x) \in X \times X : x \in X\}$ in $X \times X$, then one can show that $\lambda I_X$ is closed around any point $(x, \lambda x) \in X \times X$. This implies that $\lambda I_X$ satisfies condition (A1) in Theorem 3.1 in [2]. By (6.3), $\lambda I_X$ satisfies condition (A3) in Theorem 3.1 in [2]. Hence, $\lambda I_X$ and $G$ together satisfy all conditions in Theorem 3.1 in [2], which implies this theorem. □

In particular, if $\lambda = 1$ in Theorem 6.2, then, we immediately obtain an existence theorem of stochastic fixed-point variable for the considered set-valued mapping $G: X \times S \rightrightarrows X$.

**Corollary 6.3.** *Let $X$ be an Asplund Banach space and let $(S, \tau, \mu)$ be a topological probability space. Let $G: X \times S \rightrightarrows X$ be a set-valued mapping. Suppose that $G$ satisfies condition $(A2)'$ in Theorem 6.2 with respect to a point $(\bar{x}, \bar{s}) \in X \times S$. Then, for each $\alpha$ with $1 > \alpha > l$, there exist a neighborhood $W_\alpha \subset S$ of $\bar{s}$ and a stochastic fixed-point variable $\sigma_\alpha: W_\alpha \to X$ such that*

$$\sigma_\alpha(s) \in G(\sigma_\alpha(s), s) \quad \text{and} \quad \|\sigma_\alpha(s) - \bar{x}\| \leq \frac{\text{dist}(\bar{x}, G(\bar{x}, s))}{\alpha - l}, \quad \text{for any } s \in W_\alpha.$$

*Proof.* Take $\lambda = 1$ in Theorem 6.2, then, this corollary is proved immediately. □

Next, we consider some special cases of the set-valued mapping $G: X \times S \rightrightarrows X$ appeared in Theorem 6.2 and Corollary 6.3. Suppose that $G: X \times S \rightrightarrows X$ is defined by

$$G(x, s) = H(x) + w(s), \text{ for any } (x, s) \in X \times S, \tag{6.4}$$

where, $H: X \rightrightarrows X$ is a set valued mapping and $w: S \to X$ is a single-valued mapping, which could be considered as a generalized random variable with values in $X$. Meanwhile, $w$ can be treated as a noise in the mapping $G$. Next, we prove another existence theorem of stochastic fixed-point variable for this special set-valued mapping $G: X \times S \rightrightarrows X$ defined in (6.4).

**Corollary 6.4.** *Let $X$ be an Asplund Banach space and let $(S, \tau, \mu)$ be a topological probability space. Let $H: X \rightrightarrows X$ be a set-valued mapping and $w: S \to X$ a single-valued mapping. Suppose that $H$ and $w$ satisfy the following conditions with respect to a point $(\bar{x}, \bar{s}) \in X \times S$*

(A2)″ *There is neighborhoods $U$ of $\bar{x}$ in $X$, and $O$ of $\bar{s} \in S$ as well as a number $l$ with $0 < l < 1$ such that the multifunction $H$ is Lipschitz-like on $U$ with the uniform modulus $l$, while the multifunction $s \to H(\bar{x}) + w(s)$ is lower/inner semicontinuous at $\bar{s}$.*

*Then, for each $\alpha$ with $1 > \alpha > l$, there exist a neighborhood $W_\alpha \subset S$ of $\bar{s}$ and a stochastic fixed-point variable $\sigma_\alpha: W_\alpha \to X$ such that*

$$\sigma_\alpha(s) \in H(\sigma_\alpha(s)) + w(s) \quad \text{and} \quad \|\sigma_\alpha(s) - \bar{x}\| \leq \frac{\text{dist}(\bar{x}, H(\bar{x}) + w(s))}{\alpha - l}, \text{ for any } s \in W_\alpha.$$

*Proof.* For the given set-valued mapping $H: X \rightrightarrows X$ and single-valued mapping $w: S \to X$ in this corollary, define a set-valued mapping $G: X \times S \rightrightarrows X$ defined in (6.4). For any modulus $l \geq 0$, and for any $s \in S$, we have

$$H(x) + w(s) \subset H(u) + w(s) + l\|x - u\|\mathbb{B}_X \iff H(x) \subset H(u) + l\|x - u\|\mathbb{B}_X, \text{ for all } x, u \in U.$$

This implies that, the multifunction $H$ is Lipschitz-like on $U$ with the uniform modulus $l$ if and only if, $G(\cdot, s)$ is Lipschitz-like on $U$ for each $s \in O$ with the uniform modulus $l$. By (6.4), the condition in $(A2)''$ that multifunction $s \to H(\bar{x}) + w(s)$ is lower/inner semicontinuous at $\bar{s}$ is equivalent to that the multifunction $s \to G(\bar{x}, s)$ is lower/inner semicontinuous at $\bar{s}$. Hence, the set-valued mapping $G: X \times S \rightrightarrows X$ satisfies all conditions in Corollary 6.3, which implies this corollary. □

Next, we consider single-valued mappings, which are considered special cases of set-valued mappings with singleton values.

**Corollary 6.5.** *Let $X$ be an Asplund Banach space and let $(S, \tau, \mu)$ be a topological probability space. Let $g: X \times S \to X$ be a single-valued mapping. Suppose that $g$ satisfies the following conditions with respect to a point $(\bar{x}, \bar{s}) \in X \times S$*

(a2) *There are neighborhoods $U$ of $\bar{x}$ in $X$, and $O$ of $\bar{s} \in S$ as well as a number $l$ with $0 < l < 1$ such that the mapping $g(\cdot, s)$ is Lipschitz on $U$ with the uniform modulus $l$,*

$$\|g(x, s) - g(u, s)\| \leq l\|x - u\|, \text{ for all } x, u \in U, \text{ and for all } s \in O,$$

*while the mapping $s \to g(\bar{x}, s)$ is lower semicontinuous at $\bar{s}$.*

*Then, for any $\lambda, \alpha$ with $1 \geq \lambda > \alpha > l$, there exist a neighborhood $W_{\lambda\alpha} \subset S$ of $\bar{s}$ and a single-valued mapping $\sigma_{\lambda\alpha}: W_{\lambda\alpha} \to X$ such that*

$$\lambda\sigma_{\lambda\alpha}(s) = g(\sigma_{\lambda\alpha}(s), s) \quad \text{and} \quad \|\sigma_{\lambda\alpha}(s) - \bar{x}\| \leq \frac{\|\lambda\bar{x} - g(\bar{x}, s)\|}{\alpha - l}, \text{ for any } s \in W_{\lambda\alpha}.$$

*In particular, for each $\alpha$ with $1 > \alpha > l$, there exist a neighborhood $W_\alpha \subset S$ of $\bar{s}$ and a stochastic fixed-point variable $\sigma_\alpha: W_\alpha \to X$ such that*

$$\sigma_\alpha(s) = g(\sigma_\alpha(s), s) \quad \text{and} \quad \|\sigma_\alpha(s) - \bar{x}\| \leq \frac{\|\bar{x} - g(\bar{x}, s)\|}{\alpha - l}, \text{ for any } s \in W_\alpha.$$

**Observations 6.6.** Let $X$ be an Asplund Banach space and let $(S, \tau, \mu)$ be a topological probability space. For a given set-valued mapping $G: X \times S \rightrightarrows X$, there are several important questions for consideration, which are regarding to stochastic fixed-point problems associated with $G$ at a point $(\bar{x}, \bar{s}) \in X \times S$.

(I) For a given $\alpha$ with $1 > \alpha > l$, does there exist a neighborhood $W_\alpha \subset S$ of $\bar{s}$ and a stochastic fixed-point variable $\sigma_\alpha: W_\alpha \to X$ satisfying (6.1)

$$\sigma_\alpha(s) \in G(\sigma_\alpha(s), s), \text{ for any } s \in W_\alpha?$$

(II) If the answer for (I) is yes, then, does the stochastic fixed-point variable $\sigma_\alpha: W_\alpha \to X$ satisfy the estimation in (6.2)
$$\|\sigma_\alpha(s) - \bar{x}\| \leq \frac{\text{dist}(\bar{x}, G(\bar{x}, s))}{\alpha - l}, \text{ for any } s \in W_\alpha?$$

(III) For another given $\beta$ with $1 > \beta > l$ and $\beta \neq \alpha$, if there exist a neighborhood $W_\beta \subseteq S$ of $\bar{s}$ and a stochastic fixed-point variable $\sigma_\beta: W_\beta \to X$ such that
$$\sigma_\beta(s) \in G(\sigma_\beta(s), s), \text{ for any } s \in W_\beta,$$
then, $W_\beta = W_\alpha$?

We provide some examples to consider the above issues regarding to stochastic fixed-point problems. We start at single-valued mappings.

**Example 6.7.** Let $X = (-\infty, \infty)$ and $S = [0, 1]$, which are equipped with the standard topology. Let the probability measure on $S = [0, 1]$ equal the standard Lebesgue measure on $[0, 1]$. Define a single-valued mapping $g: X \times S \to X$ by
$$g(x, s) = \frac{1}{4}x^2 + s, \text{ for any } (x, s) \in X \times S.$$

Consider a point $(0, \bar{s}) \in X \times S$, with an arbitrary given $\bar{s} \in S$. Let $U = (-1, 1)$, which is a neighborhood of $0$ in $X$. Let $O = S$ as a special neighborhood of $\bar{s}$ in $S$. For any $s \in S$, we have
$$|g(x, s) - g(u, s)| \leq \frac{1}{2}|x - u|, \text{ for any } x, u \in U.$$

Hence, $g$ is Lipschitz-like on $U$ with the uniform modulus $l = \frac{1}{2}$, while the mapping $s \to g(0, s)$ is continuous at $\bar{s}$ (it is continuous on $S$). The mapping $g: X \times S \to X$ satisfies all conditions in Corollary 6.5. Then,

(i) For any $\alpha$ with $1 > \alpha > \frac{1}{2}$, $g$ possess a stochastic fixed-point variable on $W_\alpha = S$ with
$$\sigma_\alpha(s) = 2(1 - \sqrt{1-s}) \quad \text{and} \quad \sigma_\alpha(s) = g(\sigma_\alpha(s), s), \text{ for any } s \in S.$$

At the point $(0, \bar{s}) \in X \times S$, by $l = \frac{1}{2}$, for any $\alpha$ with $1 > \alpha > \frac{1}{2}$, we estimate
$$|\sigma_\alpha(s) - 0| = 2(1 - \sqrt{1-s}) \leq 2s = \frac{|0 - g(0,s)|}{1 - \frac{1}{2}} \leq \frac{|0 - g(0,s)|}{\alpha - l}, \text{ for any } s \in S.$$

This implies that, for any $\alpha$ with $1 > \alpha > \frac{1}{2}$, the stochastic fixed-point variable $\sigma_\alpha$ satisfies the inequality (6.2) with the neighborhood $W_\alpha = S$ of $\bar{s}$,
$$|\sigma_\alpha(s) - 0| \leq \frac{|0 - g(0,s)|}{\alpha - l}, \text{ for any } s \in W_\alpha = S.$$

(ii) At the considered point $(0, \bar{s}) \in X \times S$, with the neighborhood $U = (-1,1)$ of 0 in $X$ and the neighborhood $O = S$ of $\bar{s}$, for any $\beta$ with $1 > \beta > \frac{1}{2}$, we can show that $g$ possess another stochastic fixed-point variable $\zeta_\beta: S \to X$ such that

$$\zeta_\beta(s) = 2(1 + \sqrt{1-s}) \quad \text{and} \quad \zeta_\beta(s) = g(\zeta_\beta(s), s), \text{ for any } s \in S.$$

However, by $|0 - g(0,s)| = s$, we have

$$|\zeta_\beta(s) - 0| = 2(1 + \sqrt{1-s}) \nleq \frac{|0-g(0,s)|}{\beta-l}, \text{ for ALL } s \in W_\alpha = S.$$

(iii) At a point $(0, \bar{s}) \in X \times S$ with $\bar{s} > 0$. Let $U = (-1,1)$, which is a neighborhood of 0 in $X$ and $O = S$ as a neighborhood of $\bar{s}$ in $S$. Take an arbitrary $\bar{t}$ with $0 < \bar{t} < \bar{s}$. Then, $W_\gamma = [\bar{t}, 1]$ is a neighborhood of $\bar{s}$ in $S$. For any $\gamma$ satisfying

$$1 > \frac{1}{2} + \frac{\bar{t}}{2(1+\sqrt{1-\bar{t}})} \geq \gamma > \frac{1}{2},$$

$g$ possess another stochastic fixed-point variable $\lambda_\gamma: [\bar{t}, 1] \to X$ with $W_\gamma = [\bar{t}, 1]$.

$$\lambda_\gamma(s) = 2(1 + \sqrt{1-s}) \quad \text{and} \quad \lambda_\gamma(s) = g(\lambda_\gamma(s), s), \text{ for any } s \in S.$$

Furthermore, the following inequality is satisfied

$$|\lambda_\gamma(s) - 0| \leq \frac{|0-g(0,s)|}{\gamma-l}, \text{ for every } s \in W_\gamma = [\bar{t}, 1]. \tag{6.5}$$

Proof of (6.5). Under the condition $\frac{1}{2} + \frac{\bar{t}}{2(1+\sqrt{1-\bar{t}})} \geq \gamma > \frac{1}{2}$, we have

$$|\lambda_\gamma(s) - 0| = 2(1 + \sqrt{1-s})$$

$$\leq 2(1 + \sqrt{1-\bar{t}}) = \frac{\bar{t}}{\frac{1}{2} + \frac{\bar{t}}{2(1+\sqrt{1-\bar{t}})} - \frac{1}{2}}$$

$$\leq \frac{\bar{t}}{\gamma - \frac{1}{2}} \leq \frac{s}{\gamma - \frac{1}{2}} = \frac{|0-g(0,s)|}{\gamma-l}, \text{ for any } s \in W_\gamma = [\bar{t}, 1].$$

**Example 6.8.** Let $X$ and $S$ be given as in Example 6.7. Define a single-valued mapping $g: X \times S \to X$ by

$$g(x, s) = \frac{1}{4}x^2 s, \text{ for any } (x, s) \in X \times S.$$

Consider a point $(0, \bar{s}) \in X \times S$, with an arbitrary given $\bar{s} \in S$. Let $U = (-1,1)$ and $O = S$. For any $s \in S$, we have

$$|g(x,s) - g(u,s)| \leq \frac{1}{2}|x - u|, \text{ for any } x, u \in U.$$

The mapping $g: X \times S \to X$ satisfies all conditions in Corollary 6.5. Then,

(i) For any $\alpha$ with $1 > \alpha > \frac{1}{2}$, $g$ possess a stochastic fixed-point variable $\sigma_\alpha: S \to X$ with

$$\sigma_\alpha(s) = 0 \quad \text{and} \quad \sigma_\alpha(s) = g(\sigma_\alpha(s), s), \text{ for any } s \in S.$$

At the point $(0, \bar{s}) \in X \times S$, by $l = \frac{1}{2}$, for any $\alpha$ with $1 > \alpha > \frac{1}{2}$, the stochastic fixed-point variable $\sigma_\alpha$ satisfies the following inequality

$$|\sigma_\alpha(s) - 0| = 0 = \frac{0}{\alpha - l} = \frac{|0 - g(0,s)|}{\alpha - l}, \text{ for any } s \in S.$$

(ii) For any $\beta$ with $1 > \beta > \frac{1}{2}$, $g$ possess another stochastic fixed-point variable $\zeta_\beta: W_\beta \to X$ with $W_\beta = (0, 1]$ such that

$$\zeta_\beta(s) = \frac{4}{s} \quad \text{and} \quad \zeta_\beta(s) = g(\zeta_\beta(s), s), \text{ for any } s \in W_\beta = (0, 1].$$

However, by $|0 - g(0, s)| = 0$, we have

$$|\zeta_\beta(s) - 0| = \frac{4}{s} \not\leq 0 = \frac{|0 - g(0,s)|}{\beta - l}, \text{ for every } s \in W_\beta = (0, 1].$$

**Example 6.9.** Let $(X, \|\cdot\|) = (\mathbb{R}^2, \|\cdot\|)$ be the 2-d Euclidean space with the standard $l_2$-norm $\|\cdot\|$. Let $S = \mathbb{R} = (-\infty, \infty)$ that is equipped with the standard Boral $\sigma$-field $\mathcal{A}$, on which the probability measure $\mu$ is defined by

$$\mu(A) = \frac{1}{\sqrt{2\pi}} \int_A e^{-\frac{t^2}{2}} dt, \text{ for any } A \in \mathcal{A}.$$

For any points $a, b \in \mathbb{R}^2$, let $\overline{a, b}$ denote the closed line segment in $\mathbb{R}^2$ with ending points $a$ and $b$. Define a set-valued mapping $G: \mathbb{R}^2 \times \mathbb{R} \rightrightarrows \mathbb{R}^2$, for any $((w, v), s) \in \mathbb{R}^2 \times \mathbb{R}$, by

$$G((w,v), s) = \overline{\left(\tfrac{1}{4}w, 0\right), \left(\tfrac{1}{4}w, \tfrac{1}{4}\sqrt{1+v^2}\right)} + (s^2, |s|), \text{ for any } ((w,v), s) \in \mathbb{R}^2 \times \mathbb{R}. \quad (6.6)$$

For the point $((0,0), 0) \in \mathbb{R}^2 \times \mathbb{R}$ and the special neighborhoods $U = \mathbb{R}^2$ of $(0,0)$ in $\mathbb{R}^2$, and $O = \mathbb{R}$ of $0 \in S$, we show that the multifunction $G(\cdot, s)$ is Lipschitz-like on $\mathbb{R}^2$ (relative to $\mathbb{R}^2$) for each $s \in \mathbb{R}$ with the uniform modulus $l = \frac{1}{2}$. To this end, for $((w, v), s), ((p, q), s) \in \mathbb{R}^2 \times \mathbb{R}$, we need to prove that

$$G((w,v), s) \subset G((p,q), s) + \tfrac{1}{2}\|(w,v) - (p,q)\|\mathbb{B}_{\mathbb{R}^2}, \text{ for any } s \in \mathbb{R}. \quad (6.7)$$

By definition of in (6.6), we have

$$G((p,q), s) = \overline{\left(\tfrac{1}{4}p, 0\right), \left(\tfrac{1}{4}p, \tfrac{1}{4}\sqrt{1+q^2}\right)} + (s^2, |s|), \text{ for } ((p,q), s) \in \mathbb{R}^2 \times \mathbb{R}.$$

To prove (6.7), it is sufficient to show that, for any $\left(\tfrac{1}{4}w, \lambda\right) + (s^2, |s|) \in G((w,v), s)$, with $0 \leq \lambda \leq \tfrac{1}{4}\sqrt{1+v^2}$, there is $\beta$ with $0 \leq \beta \leq \tfrac{1}{4}\sqrt{1+q^2}$ such that, for any $s \in \mathbb{R}$,

$$\left\|\left(\tfrac{1}{4}w,\lambda\right)+(s^2,|s|)\right)-\left(\left(\tfrac{1}{4}p,\beta\right)+(s^2,|s|)\right)\right\| \leq \tfrac{1}{2}\|(w,v)-(p,q)\|. \tag{6.8}$$

The proof of (6.8) is divided into two cases with respect to $\lambda$ with $0 \leq \lambda \leq \tfrac{1}{4}\sqrt{1+v^2}$.

Case 1. $0 \leq \lambda \leq \min\{\tfrac{1}{4}\sqrt{1+v^2}, \tfrac{1}{4}\sqrt{1+q^2}\}$. In this case, we take $\beta = \lambda$, which satisfies

$$\left(\tfrac{1}{4}w,\lambda\right)+(s^2,|s|) \in G((w,v),s) \quad \text{and} \quad \left(\tfrac{1}{4}p,\lambda\right)+(s^2,|s|) \in G((p,q),s), \text{ for any } s \in \mathbb{R}.$$

For any $s \in \mathbb{R}$, we calculate

$$\left\|\left(\left(\tfrac{1}{4}w,\lambda\right)+(s^2,|s|)\right)-\left(\left(\tfrac{1}{4}p,\lambda\right)+(s^2,|s|)\right)\right\| = \tfrac{1}{4}|w-p| \leq \tfrac{1}{4}\|(w,v)-(p,q)\|.$$

This proves (6.8) for any $\lambda$ with $0 \leq \lambda \leq \min\{\tfrac{1}{4}\sqrt{1+v^2}, \tfrac{1}{4}\sqrt{1+q^2}\}$.

Case 2. (In case if there are some $\lambda$ satisfying) $\min\{\tfrac{1}{4}\sqrt{1+v^2}, \tfrac{1}{4}\sqrt{1+q^2}\} < \lambda \leq \tfrac{1}{4}\sqrt{1+v^2}$. It implies $\min\{\tfrac{1}{4}\sqrt{1+v^2}, \tfrac{1}{4}\sqrt{1+q^2}\} = \tfrac{1}{4}\sqrt{1+q^2}$. In this case, we take $\beta = \tfrac{1}{4}\sqrt{1+q^2}$. For any $\tfrac{1}{4}\sqrt{1+q^2} < \lambda \leq \tfrac{1}{4}\sqrt{1+v^2}$, by mean value theorem, for any $s \in \mathbb{R}$, one has

$$\left\|\left(\left(\tfrac{1}{4}w,\lambda\right)+(s^2,|s|)\right)-\left(\left(\tfrac{1}{4}p,\beta\right)+(s^2,|s|)\right)\right\|$$

$$= \left\|\left(\left(\tfrac{1}{4}w,\lambda\right)+(s^2,|s|)\right)-\left(\left(\tfrac{1}{4}p,\tfrac{1}{4}\sqrt{1+q^2}\right)+(s^2,|s|)\right)\right\|$$

$$= \left\|\left(\tfrac{1}{4}w,\lambda\right)-\left(\tfrac{1}{4}p,\tfrac{1}{4}\sqrt{1+q^2}\right)\right\|$$

$$\leq \left\|\left(\tfrac{1}{4}w,\tfrac{1}{4}\sqrt{1+v^2}\right)-\left(\tfrac{1}{4}p,\tfrac{1}{4}\sqrt{1+q^2}\right)\right\|$$

$$\leq \tfrac{1}{4}|w-p|+\tfrac{1}{4}\left|\sqrt{1+v^2}-\sqrt{1+q^2}\right|$$

$$= \tfrac{1}{4}|w-p|+\tfrac{1}{4}\tfrac{|t|}{\sqrt{1+t^2}}|v-q|, \text{ for some } t \text{ with } q \leq t \leq v$$

$$\leq \tfrac{1}{4}|w-p|+\tfrac{1}{4}|v-q|$$

$$\leq \tfrac{1}{2}\|(w,v)-(p,q)\|.$$

This proves (6.8) for any $\lambda$ (if any) with $\min\{\tfrac{1}{4}\sqrt{1+v^2}, \tfrac{1}{4}\sqrt{1+q^2}\} < \lambda \leq \tfrac{1}{4}\sqrt{1+v^2}$. Hence, (6.7) is proved; and therefore, the multifunction $G(\cdot, s)$ is Lipschitz-like on $\mathbb{R}^2$ (relative to $\mathbb{R}^2$) for each $s \in \mathbb{R}$ with the uniform modulus $l = \tfrac{1}{2}$. For any given point $x \in \mathbb{R}^2$, the multifunction $s \to G(x,s)$ is continuous on $\mathbb{R}$. Hence, $G$ satisfies all conditions in Corollary 6.3. which implies, for each $\alpha$ with $1 > \alpha > l\,(=\frac{1}{2})$, there exist a neighborhood $W_\alpha \subset \mathbb{R}$ of $0$ and a stochastic fixed-point variable $\sigma_\alpha\colon W_\alpha \to \mathbb{R}^2$ such that

$$\sigma_\alpha(s) \in G(\sigma_\alpha(s), s) \tag{6.9}$$

and

$$\|\sigma_\alpha(s) - (0,0)\| \le \frac{\mathrm{dist}((0,0), G((0,0),s))}{\alpha - l}, \quad \text{for any } s \in W_\alpha. \tag{6.10}$$

Actually, we precisely find infinitely many examples of stochastic fixed-point variables on the special neighborhood $\mathbb{R}$ of $0$ satisfying both (6.9) and (6.10).

For any $\alpha$ with $1 > \alpha > \frac{1}{2}$ with the special neighborhood $\mathbb{R}$ of $0$, for arbitrarily given $\lambda$ with $1 \le \lambda < \frac{4}{3}$, we define $\sigma_{(\lambda)}\colon \mathbb{R} \to \mathbb{R}^2$ by (it does not depend on $\alpha$)

$$\sigma_{(\lambda)}(s) = \left(\tfrac{4}{3}s^2, \lambda|s|\right), \quad \text{for any } s \in \mathbb{R}.$$

Note that in (6.9) and (6.10), the foot notes $\alpha$ of the stochastic fixed-point variable $\sigma_\alpha$ is related with the Lipschitz-like modulus of $G$. Meanwhile, the given $\lambda$ with $1 \le \lambda < \frac{4}{3}$ is not related to the Lipschitz-like modulus of $G$. To distinguish the difference, we use $\sigma_{(\lambda)}$ for the stochastic fixed-point variable, which does not depend on $\alpha$. Then, we have

$$\begin{aligned}
G(\sigma_{(\lambda)}(s), s) &= G\left(\left(\tfrac{4}{3}s^2, \lambda|s|\right), s\right) \\
&= \overline{\left(\tfrac{1}{4}\tfrac{4}{3}s^2, 0\right), \left(\tfrac{1}{4}\tfrac{4}{3}s^2, \tfrac{1}{4}\sqrt{1 + \lambda^2|s|^2}\right)} + (s^2, |s|) \\
&= \overline{\left(\tfrac{1}{3}s^2, 0\right), \left(\tfrac{1}{3}s^2, \tfrac{1}{4}\sqrt{1 + \lambda^2 s^2}\right)} + (s^2, |s|) \\
&= \overline{\left(\tfrac{4}{3}s^2, |s|\right), \left(\tfrac{4}{3}s^2, |s| + \tfrac{1}{4}\sqrt{1 + \lambda^2 s^2}\right)} \\
&\ni \left(\tfrac{4}{3}s^2, \lambda|s|\right) = \sigma_{(\lambda)}(s), \quad \text{for any } s \in \mathbb{R}.
\end{aligned}$$

The above inclusion is from the assumption $1 \le \lambda < \frac{4}{3}$, which implies

$$|s| \le \lambda|s| < \tfrac{4}{3}|s| < |s| + \tfrac{1}{4}\sqrt{1 + \lambda^2 s^2}, \text{ for any } s \in \mathbb{R}.$$

We obtain that $\sigma_{(\lambda)}(s)$ satisfies (6.9), for any $1 \le \lambda < \frac{4}{3}$,

$$\sigma_{(\lambda)}(s) = \left(\tfrac{4}{3}s^2, \lambda|s|\right) \in G(\sigma_{(\lambda)}(s), s), \text{ for any } s \in \mathbb{R}.$$

This implies that, for any $1 \le \lambda < \frac{4}{3}$, $\sigma_{(\lambda)}$ is a stochastic fixed-point variable on $\mathbb{R}$, which does not depend on the Lipschitz-like modulus of $G$. Before checking that the stochastic fixed-point variable $\sigma_{(\lambda)}$ satisfies (6.10), we firstly estimate $\text{dist}((0,0), G((0,0), s))$. We calculate

$$G((0,0), s) = \overline{\left(\tfrac{1}{4}, 0, 0\right), \left(\tfrac{1}{4}, 0, \tfrac{1}{4}\sqrt{1+0^2}\right)} + (s^2, |s|)$$

$$= \overline{(0,0), \left(0, \tfrac{1}{4}\right)} + (s^2, |s|)$$

$$= \overline{(s^2, |s|), \left(s^2, \tfrac{1}{4} + |s|\right)}.$$

This is a vertical closed line segment with ending points $(s^2, |s|)$ and $\left(s^2, \tfrac{1}{4} + |s|\right)$. We have

$$\text{dist}\left((0,0), G((0,0), s)\right)$$

$$= \min\left\{\|(0,0) - (s^2, \beta)\|: |s| \leq \beta \leq \tfrac{1}{4} + |s|\right\}$$

$$= \|(0,0) - (s^2, |s|)\|$$

$$= |s|\sqrt{s^2 + 1}, \text{ for any } s \in \mathbb{R}.$$

Then, for any $\alpha$ with $1 > \alpha > \tfrac{1}{2}$ and for any $\lambda$ with $1 \leq \lambda < \tfrac{4}{3}$, we estimate

$$\|\sigma_{(\lambda)}(s) - (0,0)\| = \left\|\left(\tfrac{4}{3}s^2, \lambda|s|\right)\right\|$$

$$= |s|\sqrt{\tfrac{16}{9}s^2 + \lambda^2} < \tfrac{4}{3}|s|\sqrt{s^2 + 1}$$

$$< 2|s|\sqrt{s^2 + 1}$$

$$= \frac{|s|\sqrt{s^2+1}}{1 - \tfrac{1}{2}} < \frac{|s|\sqrt{s^2+1}}{\alpha - \tfrac{1}{2}}$$

$$= \frac{\text{dist}((0,0), G((0,0), s))}{\alpha - \tfrac{1}{2}}, \text{ for any } s \in \mathbb{R}.$$

This proves that the stochastic fixed-point variable $\sigma_{(\lambda)}$ satisfies (6.10), for any $\alpha$ with $1 > \alpha > \tfrac{1}{2}$ and for any $\lambda$ with $1 \leq \lambda < \tfrac{4}{3}$, which ensures that the set-valued mapping $G: \mathbb{R}^2 \times \mathbb{R} \rightrightarrows \mathbb{R}^2$ has infinitely many stochastic fixed-point variables $\sigma_{(\lambda)}$ on $\mathbb{R}$. All the stochastic fixed-point variables $\sigma_{(\lambda)}$ do not depend on the Lipschitz-like modulus of $G$.

**Acknowledgments** The author is very grateful to Professor Boris S. Mordukhovich for his kind communications, valuable suggestions and enthusiasm encouragements in the development stage of this paper. The author sincerely thanks Professor Robert Mendris for his valuable suggestions, which improved the presentation of this paper.